\documentclass{amsart}
\usepackage{amsthm,amssymb}

\newtheorem{proposition}{Proposition}[section]
\newtheorem{theorem}[proposition]{Theorem}
\newtheorem{corollary}[proposition]{Corollary}
\newtheorem{lemma}[proposition]{Lemma}
\newtheorem{conjecture}[proposition]{Conjecture}

\theoremstyle{definition}
\newtheorem{remark}[proposition]{Remark}
\newtheorem*{definition}{Definition}

\newtheorem*{notation}{Notation}
\newtheorem{example}[proposition]{Example}
\newtheorem{algorithm}[proposition]{Algorithm}

\newcounter{alphcount}
\renewcommand{\thealphcount}{(\alph{alphcount})}
\newenvironment{alph-list}{\begin{list}{\thealphcount}{\usecounter{alphcount}}}%
                        {\end{list}} 

\newcounter{romcount}
\renewcommand{\theromcount}{(\roman{romcount})}
\newenvironment{rom-list}{\begin{list}{\theromcount}{\usecounter{romcount}}}%
                        {\end{list}} 

\newcounter{caseno}
\renewcommand{\thecaseno}{\arabic{caseno}}

\newcounter{stepno}
\renewcommand{\thestepno}{\arabic{stepno}}
\newenvironment{steplist}{\begin{list}{\em Step \thestepno :}
				   {\usecounter{stepno}}%
		      }{\end{list}}

\newcommand{\bdmap}{\partial}

\newcommand{\cC}{\mathcal C}
\newcommand{\cI}{\mathcal I}
\newcommand{\cB}{\mathcal B}
\newcommand{\cF}{\mathcal F}
\newcommand{\cV}{\mathcal V}

\newcommand{\calS}{\mathcal S}
\newcommand{\Sf}{\mathcal D}

\newcommand{\K}{\mathcal K}
\newcommand{\IN}{IN}
\newcommand{\s}{\mathbf{s}}
\newcommand{\dd}{\mathbf{d}}

\newcommand{\clM}[1]{\overline{#1}}
\newcommand{\ostar}{\circ}
\newcommand{\ppi}[1]{\pi_{#1}}
\newcommand{\pibar}[1]{\overline{\pi}_{#1}}
\newcommand{\chit}[1]{\tilde{\chi}(#1)}
\newcommand{\chitcx}[1]{\tilde{\chi}(#1)}
\newcommand{\chitempty}{\tilde{\chi}}
\newcommand{\abs}[1]{\lvert#1\rvert}
\newcommand{\setm}{\backslash}

\newcommand{\cx}[1]{\Delta(#1)}
\newcommand{\coeff}[1]{[#1]}
\newcommand{\pskel}[1]{^{[#1]}}
\newcommand{\skel}[1]{^{(#1)}}
\newcommand{\majby}{\trianglelefteq}
\newcommand{\id}{\mathrm{id}}

\newcommand{\reals}{\mathbb R}
\newcommand{\disun}{\mathbin{\dot{\cup}}} 
\newcommand{\condns}[2]{\substack{#1 \\ #2}}

\newcommand{\blambda}{\boldsymbol{\lambda}}
\newcommand{\bmu}{\boldsymbol{\mu}}

\DeclareMathOperator{\Spec}{Spec}
\DeclareMathOperator{\rk}{rk}
\DeclareMathOperator{\cl}{cl}
\DeclareMathOperator{\ci}{ci}
\DeclareMathOperator{\bo}{bo}

\DeclareMathOperator{\lk}{lk}
\DeclareMathOperator{\st}{st}
\DeclareMathOperator{\bd}{bd}

\newcommand{\ie}{{\em i.e.}}
\newcommand{\eg}{{\em e.g.}}
\newcommand{\cf}{{\em cf.}}
\newcommand{\Cf}{{\em Cf.}}

\begin{document}

\title[Laplacians of matroids and shifted complexes]{A common recursion for Laplacians of
matroids and shifted simplicial complexes}

\author{Art M. Duval}
\email{artduval@math.utep.edu}
\address{Department of Mathematical Sciences\\
         University of Texas at El Paso\\
         El Paso, TX 79968-0514}



\begin{abstract}
A recursion due to Kook expresses the Laplacian eigenvalues of a
matroid $M$ in terms of the eigenvalues of its deletion $M-e$
and contraction $M/e$ by a fixed element $e$, and an error term.  We
show that this error term is given simply by the Laplacian eigenvalues
of the pair $(M-e, M/e)$.  
We further show that by suitably generalizing deletion and contraction
to arbitrary simplicial complexes, the Laplacian eigenvalues of
shifted simplicial complexes satisfy this exact same recursion.

We show that the class of simplicial complexes satisfying this
recursion is closed under a wide variety of natural operations, and
that several specializations of this recursion reduce to basic
recursions for natural invariants.

We also find a simple formula for the Laplacian eigenvalues of an
arbitrary pair of shifted complexes in terms of a kind of generalized
degree sequence.

\end{abstract}

\maketitle


\section{Introduction}\label{se:intro}

The independence complex of matroids and shifted simplicial complexes
are two of only four types of simplicial complexes whose combinatorial
Laplacians $L=\partial\partial^{*}+\partial^{*}\partial$ are known to
have only integer eigenvalues (see Kook, Reiner, and Stanton
\cite{KookReinerStanton}, and \cite{DuvalReiner}, respectively).  The
other two types, which will not concern us further, are matching
complexes of complete graphs \cite{DongWachs} and chessboard complexes
\cite{FriedmanHanlon}.  More information and background about the
combinatorial Laplacian and its eigenvalues may be found in Section
\ref{se:laplace} and \cite{DuvalReiner,Friedman,KookReinerStanton}.
Our main result (Theorems \ref{th:matroids} and \ref{th:shifted}) is
another, more striking, similarity between the Laplacian eigenvalues
of matroids and shifted complexes: they satisfy the exact same
recursion, which we call the {\em spectral recursion}, equation
\eqref{eq:big}. This recursion is stated in terms of the {\em spectrum
polynomial}, a natural generating function for Laplacian eigenvalues,
defined in equation \eqref{eq:S.defn}.

The Tutte polynomial $T_{M}$ of a matroid $M$ satisfies the recursion
$T_{M} = T_{M-e}+T_{M/e}$, when $e$ is neither a loop nor an isthmus,
and where $M-e$ and $M/e$ denote the deletion and contraction,
respectively, of $M$ with respect to ground element $e$.  When Kook,
Reiner, and Stanton proved that the Laplacian spectrum of a matroid is
integral, they also speculated on the existence of a Tutte
polynomial-like recursion for the spectrum polynomial of a matroid
$M$, though possibly with a third ``error'' term, besides the deletion
and contraction, on the right-hand side \cite[Question
3]{KookReinerStanton}.  Kook \cite{Kook} found such a recursion, but
the error term in his formulation is somewhat complicated to state,
with two cases depending on whether or not the ground element $e$ is a
closed element in $M$.  Subsequently, Kook and Reiner (private
communication) asked if this error term might be just the spectrum
polynomial of the matroid pair $(M-e,M/e)$.

One of our main results (Theorem \ref{th:matroids}) is that Kook and
Reiner's conjecture is true, that is, the spectrum polynomial of $M$
can be expressed simply in terms of the spectrum polynomials of $M-e$,
$M/e$, and $(M-e,M/e)$.  This is the spectral recursion.  We show,
furthermore, by suitably generalizing the definitions of deletion and
contraction from matroids to arbitrary simplicial complexes (Section
\ref{se:laplace}), that shifted complexes also satisfy the spectral
recursion (Theorem \ref{th:shifted}).

This raises the natural question: What is the largest class of
simplicial complexes, necessarily a common generalization of matroids
and shifted complexes, satisfying the spectral recursion?  We will see
that this class is closed under the operations of join, skeleta,
Alexander dual, and disjoint union (Corollaries
\ref{th:join.ok.all}, \ref{th:skeleta.preserve.all},
\ref{th:spectral.Alexander.all}, and \ref{th:union.preserve.all},
respectively).  We might hope that it is closed also under deletion
and contraction, as matroids and shifted complexes each are.  In the
same vein, it may be worthwhile to restrict our attention to those
complexes that are also Laplacian integral.  Unfortunately, no hint to
determining this common generalization is apparent in the proofs of
either Laplacian integrality or the spectral recursion, which are each
rather different for matroids and shifted complexes.

Jarrah and Laubenbacher \cite{JarrahLaubenbacher} examined another
property shared by matroids and shifted complexes.  Klivans
\cite{Klivans} has characterized simplicial complexes that are
simultaneously shifted {\em and} the matroid complex of some matroid;
this is, in some sense, the reverse of finding a natural common
generalization of matroids and shifted complexes.

The common generalization includes neither of the other known types of
Laplacian integral simplicial complexes.  Direct computations show
that the matching complex of the complete graph on 5 vertices and the
$2 \times 3$ chessboard complex both fail to satisfy the spectral recursion
with respect to any vertex.  Also excluded is the 3-edge path
(Example \ref{ex:universal}), which rules out as the common generalization
such otherwise likely candidates as vertex-decomposable
\cite{ProvanBillera}\cite[Section 11]{BjornerWachs2} or shellable
complexes \cite{BjornerWachs1,BjornerWachs2}.

A key piece of the proof that matroids satisfy the spectral recursion
is a decomposition of the Laplacian of $(M-e,M/e)$ into a direct sum of Laplacians of $M/C$'s for
all circuits $C$ containing $e$ (Lemma \ref{th:matroid.partition}).
We may combine this with the spectral recursion to express the
spectrum polynomial of a matroid completely in terms of spectrum
polynomials of smaller matroids (with no matroid pairs), which permits
a truly recursive way of computing Laplacian eigenvalues for matroids
(Remark \ref{rm:truly.recursive}).

Unfortunately, we are unable to state any formula for the Laplacian
eigenvalues of an {\em arbitrary} matroid pair (\ie, besides
$(M-e,M/e)$).  We are able, however, to use tools developed in the
proof of the spectral recursion for shifted complexes to find a simple
formula for the Laplacian eigenvalues of an arbitrary shifted
simplicial pair (Theorem \ref{th:sdT}).  This naturally generalizes a
formula for a single shifted complex \cite{DuvalReiner}; the graph
case goes back to Merris \cite{Merris}.  Similarly, we generalize a
related conjectured inequality on the Laplacian spectrum of an
arbitrary simplicial complex \cite{DuvalReiner} to an arbitrary
simplicial pair (Conjecture \ref{th:GM}); the graph case was
conjectured by Grone and Merris \cite{GroneMerris}.  Passing from
graphs to simplicial complexes in \cite{DuvalReiner} required
generalizing the well-known notion of degree sequences for graphs.
Now passing to simplicial {\em pairs}, we introduce a less than
obvious, but perfectly natural, further generalization of degree
sequence (Subsection \ref{su:degree}).

The Tutte polynomial is arguably the most important invariant of
matroid theory (see, \eg, \cite{BrylawskiOxley}).  The spectrum
polynomial shares several nice features with the Tutte polynomial,
such as being well-behaved under join (Corollary \ref{th:S.join}),
disjoint union (Lemma \ref{th:S.union}), and several dual operators
(equations \eqref{eq:S.tn.dual} and \eqref{eq:S.complement}).
Furthermore, specializations obtained by plugging in particular values
for one or the other of the variables of the spectrum polynomial
reduce it to well-known invariants.  Consequently (and now going
beyond matroids and the Tutte polynomial), in each of these
specializations, the spectral recursion holds for all simplicial
complexes $\Delta$ (not just matroids and shifted complexes), because
it reduces to a basic recursion expressing the relevant invariant for
$\Delta$ in terms of that invariant for $\Delta-e$ and $\Delta/e$
(Theorem \ref{th:homology} and Corollary \ref{th:t-1}).

In contrast to the Tutte polynomial recursion, the spectral recursion does
not need to exclude loops and isthmuses as special cases.  Indeed,
the spectral recursion holds for all complexes (not just matroids and
shifted complexes) when $e$ is a loop (Proposition \ref{th:loop.ok}) or an isthmus
(Proposition \ref{th:v.ok} and Theorem \ref{th:homology}).

Section \ref{se:laplace} contains more information about Laplacians
and the spectral recursion, including some special cases.  Sections
\ref{se:matroids} and \ref{se:shifted} are devoted to the proofs that
matroids and shifted complexes, respectively, satisfy the spectral
recursion.  The formula for eigenvalues of arbitrary shifted
simplicial pairs is developed in Section \ref{se:arb.shifted}.
Finally, in Section \ref{se:preserve}, we show that disjoint union and
several duality operators, including Alexander duality, all preserve
the property of satisfying the spectral recursion.

\section{Laplacians of simplicial pairs}\label{se:laplace}

For further background on simplicial complexes, their boundary maps
and homology groups, see, \eg, \cite[Chapter 1]{Munkres}.
If $\Delta$ and $\Delta'$ are simplicial complexes on the same ground
set of vertices, then we will say $(\Delta,\Delta')$ is a {\em
simplicial pair}, but we set $(\Delta,\Delta')=(\Gamma,\Gamma')$ when
the set differences $\Delta \setm \Delta'$ and $\Gamma \setm \Gamma'$
are equal as subsets of the power set of the ground set of vertices
(here $A \setm B$ denotes the set difference $\{a \in A: a \not\in
B\}$ between sets $A$ and $B$); more formally, then, a simplicial pair
is an equivalence class on ordered pairs of simplicial complexes.  In
all cases, definitions applying to a simplicial pair
($\Delta,\Delta'$) may be specialized to a single simplicial complex
$\Delta$, by letting $\Delta'=\emptyset$, the empty simplicial
complex.  

As usual, let $C_i = C_i(\Delta,\Delta';\reals) :=
C_i(\Delta;\mathbb{R})/C_i(\Delta';\reals)$ denote the $i$-dimensional
oriented $\reals$-chains of $(\Delta,\Delta')$, \ie, the formal
$\reals$-linear sums of oriented $i$-dimensional faces 
$[F]$ such that $F \in \Delta_i \setm \Delta'_i$,
where $\Delta_i$ denotes the set of $i$-dimensional faces of $\Delta$.
Let
$\bdmap_{(\Delta,\Delta');i}=\bdmap_i\colon C_i \rightarrow C_{i-1}$
denote the usual (signed) {\em boundary operator}.
Via the natural bases $\Delta_i \setm \Delta'_i$ and $\Delta_{i-1}
\setm \Delta'_{i-1}$ for $C_i(\Delta,\Delta';\reals)$ and
$C_{i-1}(\Delta,\Delta';\reals)$, respectively, the boundary map
$\bdmap_i$ has an adjoint map $\bdmap^*_i\colon
C_{i-1}(\Delta,\Delta'; \reals) \rightarrow C_{i}(\Delta,\Delta';
\reals)$; \ie,
the matrices representing $\bdmap$ and $\bdmap^*$ in the natural bases
are transposes of one another.

\begin{definition}
Let $L'_i=\bdmap_{i+1}\bdmap_{i+1}^*$ and $L''_i=\bdmap_{i}^*\bdmap_{i}$.
Then the {\em {\rm (}$i$-dimensional\,{\rm )} Laplacian} of
$(\Delta,\Delta')$ is the map $L_i(\Delta,\Delta') \colon
C_i(\Delta,\Delta';\reals) \rightarrow
C_i(\Delta,\Delta';\reals)$ defined by
$$
L_i = L_i(\Delta,\Delta') := L'_i + L''_i = \bdmap_{i+1}\bdmap_{i+1}^* + \bdmap_{i}^*\bdmap_{i}.
$$
\end{definition}
For more information, see, \eg, \cite{DuvalReiner,Friedman,KookReinerStanton}.
Laplacians of pairs of graphs were considered in \cite{ChungEllis}.
Each of $L'_i$ and $L''_i$ is positive semidefinite, since each is the
composition of a linear map and its adjoint.  Therefore, their sum
$L_i$ is also positive semidefinite, and so has only non-negative real
eigenvalues.  (See also Proposition \ref{th:s.s''} and
\cite[Proposition 2.1]{Friedman}.)  These eigenvalues do not depend on
the arbitrary ordering of the vertices of $\Delta$, and are thus
invariants of $(\Delta,\Delta')$; see, \eg, \cite[Remark
3.2]{DuvalReiner}.  Define $\s_i(\Delta,\Delta')$ to be the multiset
of eigenvalues of $L_i(\Delta,\Delta')$, and define
$m_\lambda(L_i(\Delta,\Delta'))$ to be the multiplicity of $\lambda$
in $\s_i(\Delta,\Delta')$.  
The single complex case ($\Delta' = \emptyset$) of the following
proposition is the first result of combinatorial Hodge theory, which
goes back to Eckmann \cite{Eckmann}.
\begin{proposition}\label{th:beta.m0}
The multiplicity of 0 as an eigenvalue of the $i$-dimensional
Laplacian $L_i$ of $(\Delta,\Delta')$ is the $i$th reduced Betti
number of $(\Delta,\Delta')$, \ie,
$$m_0(L_i(\Delta,\Delta'))=\tilde{\beta}_i(\Delta,\Delta') =
\dim_{\reals}\tilde{H}_i(\Delta,\Delta';\reals).$$
\end{proposition}
\begin{proof}
A nice summary is given in the proof
of \cite[Proposition 2.1]{Friedman}.  The usual setup is for just a
single simplicial complex (\ie, the special case $\Delta'=\emptyset$),
but only depends on the $C_i$'s and $\bdmap_i$'s forming a chain
complex ($\bdmap^2=0$), which still holds even when $\Delta' \neq \emptyset$.
(\Cf\ Proposition \ref{th:s.s''}.)
\end{proof}

A natural generating function for the Laplacian eigenvalues of a
simplicial pair $(\Delta,\Delta')$ is
\begin{equation}\label{eq:S.defn}
S_{(\Delta,\Delta')}(t,q) := \sum_{i \geq 0} t^i \sum_{\lambda \in \s_{i-1}(\Delta,\Delta')} q^\lambda
        = \sum_{i, \lambda} m_\lambda(L_{i-1}(\Delta,\Delta')) t^i q^\lambda.
\end{equation}
We call $S_{(\Delta,\Delta')}$ the {\em spectrum polynomial} of
$(\Delta,\Delta')$.  Although $S_{(\Delta,\Delta')}$ is defined for any
simplicial pair $(\Delta,\Delta')$, it is only truly a polynomial when
the Laplacian eigenvalues are not only non-negative, but integral as
well.  This will be true for the cases we are concerned with,
primarily matroids \cite{KookReinerStanton}, shifted complexes
\cite{DuvalReiner}, and shifted simplicial pairs (Theorem \ref{th:sdT}
and Remark \ref{rm:fam.to.cx}).  For the special case of a matroid, a
``spectrum polynomial'' $\Spec$ was defined, differently, in
\cite{KookReinerStanton}, but we will see later that the two
definitions agree in this case up to simple changes in indexing (see
Lemma \ref{th:KRS-step} and \cite[Corollary 18]{KookReinerStanton}).
Letting $\lambda \in \s_{i-1}$ instead of $\lambda \in \s_{i}$
simplifies the statement of some later results, notably Corollary
\ref{th:S.join}.

Recall (\eg, \cite[Section 7.3]{Bjorner}) the {\em independence
complex} $\IN(M)$ of a matroid $M$ on ground set $E$ is the simplicial
complex whose faces are the independent sets of $M$ and whose vertex
set is $E$.  (For background about matroids, see, \eg, \cite{Oxley,
Welsh, White}.)  We will sometimes use $M$ and $\IN(M)$
interchangeably, so, for instance,
$L_i(M):=L_i(\IN(M))=L_i(\IN(M),\emptyset)$ and $S_M := S_{\IN(M)} =
S_{(\IN(M),\emptyset)}$.  Similarly, if $N$ is another matroid on the
same ground set such that $\IN(N) \subseteq \IN(M)$ (\ie, $N \leq M$
in the weak order on matroids), then $L_i(M,N) = L_i(\IN(M),\IN(N))$
and $S_{(M,N)} = S_{(\IN(M),\IN(N))}$.  In this case, we say $(M,N)$
is a {\em matroid pair}.

We now naturally generalize the notion of deletion and contraction for
matroids (see \eg, \cite{Brylawski}) to arbitrary simplicial
complexes.
\begin{definition}
Let $\Delta$ be a simplicial complex on vertex set $V$, and $e \in V$. Then
the {\em deletion} of $\Delta$ with respect to $e$ is the simplicial complex
$$
\Delta - e =\{F \in \Delta\colon e \not\in F\}
$$
on vertex set $V-e$,
and the {\em contraction} of $\Delta$ with respect to $e$ is the simplicial complex
$$
\Delta/e =\{F-e\colon F \in \Delta,\ e \in F\}
$$
on vertex set $V-e$.  Note that $\Delta/e = \lk_\Delta e$, the usual
simplicial complex link \cite[Section 2]{Munkres}; we use the term
``contraction'' to highlight similarities to matroid theory.
\end{definition}

It is easy to verify that $\IN(M-e) = \IN(M)-e$ as long as $e$ is not
an isthmus of $M$, and that $\IN(M/e) = \IN(M)/e$ as long as $e$ is
not a loop of $M$.  
There is thus no confusion in 
the notational shortcuts
$S_{M-e}:=S_{\IN(M-e)}=S_{\IN(M)-e}$ and
$S_{M/e}:=S_{\IN(M/e)}=S_{\IN(M)/e}$ as long as $e$ is not an
isthmus or a loop, respectively.

Since $e$ is an isthmus of $M$ precisely when $e$ is a vertex of every
facet of $\IN(M)$, define $e$ to be an {\em isthmus} of a simplicial
complex $\Delta$ if $e$ is a vertex of every facet of $\Delta$ (so
$\Delta$ is a cone with apex $e$ -- see Subsection \ref{su:joins}).  Similarly,
since $e$ is a loop of $M$ precisely when $e$ is not a vertex of any
face of $\IN(M)$, define $e$ to be a {\em loop} of a simplicial complex
$\Delta$ if 
$e$ is in the vertex set of $\Delta$, but in no face of
$\Delta$ (even the singleton $\{e\}$ is not a face, contrary to usual simplicial complex conventions).  

Our definitions mean that if $e$ is an isthmus of simplicial complex
$\Delta$, then the deletion $\Delta - e$ equals $\Delta/e$.  (When $e$
is an isthmus of a matroid $M$, the {\em matroid} deletion $M-e$ is
left undefined in \eg, Brylawski \cite{Brylawski}, though $M-e = M/e$
in Welsh \cite[Section 4.2]{Welsh} and Oxley \cite[Corollary
3.1.25]{Oxley}.)  If $e$ is a loop of simplicial complex $\Delta$,
then the contraction $\Delta/e$ is $\emptyset$, the empty simplicial
complex.  (When $e$ is a loop of a matroid $M$, the {\em matroid}
contraction $M/e$ equals $M-e$.)

\begin{definition}
We will say that a simplicial complex $\Delta$ {\em satisfies the spectral
recursion with respect to $e$} if $e$ is a vertex of $\Delta$ and
\begin{equation}\label{eq:big}
S_\Delta(t,q) = qS_{\Delta-e}(t,q) + qtS_{\Delta/e}(t,q) + (1-q)S_{(\Delta-e,\Delta/e)}(t,q).
\end{equation}
We will say $\Delta$ {\em satisfies
the spectral recursion} if $\Delta$ satisfies the spectral recursion with respect
to every vertex in its vertex set.  (Note that Proposition \ref{th:loop.ok}
below means we need not be too particular about the vertex set of
$\Delta$.)
\end{definition}

Our main result is that
$\Delta$ satisfies the spectral recursion
when $\Delta$ is either the independence complex of a matroid (Theorem
\ref{th:matroids}) or a shifted simplicial complex (Theorem
\ref{th:shifted}), and $e$ is any vertex of $\Delta$.
We illustrate now a few special cases of the spectral recursion, which
are easy to verify, and some of which are used in later sections.

\begin{proposition}\label{th:v.ok}
The simplicial complex whose sole facet is a single vertex satisfies
the spectral recursion.
\end{proposition}

\begin{proposition}\label{th:loop.ok}
If $e$ is a loop of simplicial complex $\Delta$, then $\Delta$
satisfies the spectral recursion with respect to $e$.
\end{proposition}

Proposition \ref{th:v.ok} and Theorem \ref{th:join.ok.e} will show
that, if $e$ is an isthmus of $\Delta$, then $\Delta$ satisfies the spectral
recursion with respect to $e$.  

\begin{theorem}\label{th:homology}
If $\Delta$ is any simplicial complex, and $e$ is any vertex of
$\Delta$, then the spectral recursion
holds when $q=0$, $q=1$, $t=0$, or $t=-1$.
\end{theorem}
\begin{proof}
Plugging $q=0$ into $S$ immediately yields
$
S_{(\Delta,\Delta')}(t,0) = \sum_i t^i \tilde{\beta}_{i-1}(\Delta,\Delta')
$,
by Proposition \ref{th:beta.m0}.  Proving the spectral recursion in this
case then reduces to showing
\begin{equation}\label{eq:big.prime}
\tilde{\beta}_{i-1}(\Delta) = \tilde{\beta}_{i-1}(\Delta-e,\Delta/e),
\end{equation}
for all $i$.  This, in turn, is a consequence of the basic topology
facts 
$\tilde{\beta}_{i-1}(\Delta) = \tilde{\beta}_{i-1}(\Delta,\st_{\Delta} e)$ 
and
$(\Delta,\st_{\Delta}e) = (\Delta-e,\Delta/e)$, where $\st_{\Delta} e$
denotes the usual star of $e$ in $\Delta$, the simplicial complex
whose facets are the facets of $\Delta$ containing $e$.

Setting $q=1$, we see 
$S_{(\Delta,\Delta')}(t,1) = \sum_i (f_{i-1}(\Delta)-f_{i-1}(\Delta')) t^i$, 
where $f_i$ is the number of $i$-dimensional faces of $\Delta$, since
there are as many eigenvalues of $L_{i-1}(\Delta,\Delta')$ as there
are faces in $\Delta_{i-1} \setm \Delta'_{i-1}$ (assuming $\Delta'
\subseteq \Delta$).
It is then an easy exercise to verify that, when $q=1$, the $t^{i+1}$
coefficient of the spectral recursion reduces to the easy observation
\begin{equation}\label{eq:fTG}
f_i(\Delta) = f_i(\Delta-e) + f_{i-1}(\Delta/e).
\end{equation}

If we set $t=0$, it is easy to see that $S_\Delta(0,q) =
q^{v(\Delta)}$, where $v(\Delta)$ denotes the number of non-loop
vertices of $\Delta$.  
The spectral recursion in this case reduces to the trivial
observation that $v(\Delta)=1+v(\Delta-e)$ if $e$ is not a loop, but
$v(\Delta) = v(\Delta-e)$ if $e$ is a loop.

We will also see in Corollary \ref{th:t-1} that, when $t=-1$, the
spectral recursion reduces to an easy identity about Euler
characteristic.
\end{proof}

In the special case where $\Delta$ is a near-cone (see Subsection
\ref{su:shift.finish}) and $e$ is its apex, it is not hard to verify
that the $t^{\dim \Delta + 1}$ coefficient of the spectral recursion
reduces to \cite[Lemma 5.3]{DuvalReiner}.

The following complex is the simplest and smallest counterexample to
both Laplacian integrality and the spectral recursion.
\begin{example}\label{ex:universal}    
Let $\Delta$ be the 1-dimensional simplicial complex with vertices
$a,b,c,d$ and facets (maximal faces) $\{a,b\}$,$\{b,c\}$, and
$\{c,d\}$.  It is easy to check directly that $\Delta-e$, $\Delta/e$,
and $(\Delta-e,\Delta/e)$ are all Laplacian integral for any choice of
$e$, while $\Delta$ is not integral.  It then follows immediately that
$\Delta$ does not satisfy the spectral recursion for any choice of $e$.
\end{example}

\section{Matroids}\label{se:matroids}

In this section, we show that the independence complex of a matroid
satisfies the spectral recursion, equation \eqref{eq:big}.  The key
step of the section is a simple trick in Subsection \ref{su:partition}
to reduce the problem of computing $S_{(M-e,M/e)}$ to computing
$S_{M/C}$ for all circuits $C$ containing $e$.  Subsection
\ref{su:KRS} shows how an algorithm due to Kook, Reiner, and Stanton
\cite{KookReinerStanton} allows us to compute the spectrum polynomial
of a matroid from its combinatorial information; we also compare what
this algorithm computes for $M$, $M-e$, $M/e$, and $M/C$.  The final
steps of the calculation, which largely consist of translating to
generating functions the results of the previous subsections, are in
Subsection \ref{su:mat.finish}.

We first set our notation for matroids; for further background, and
any terms not defined here, see \cite{White}.  Let $M=M(E)$ be a
matroid on ground set $E$.  We will let $\cB=\cB(M)$, $\cI=\cI(M)$,
$\cC=\cC(M)$, and $\cF=\cF(M)$ denote the sets of bases, independent
sets, circuits, and flats (closed sets) of $M$, respectively.  If $A
\subseteq E$, let $\rk_M(A)=\rk(A)$ denote the {\em rank} of A (with
respect to $M$), and let $\clM{A} = \cl_{M}(A)$ denote the {\em
closure} of $A$ (with respect to $M$).  We will often write $V$ for
$M(V)$ in the special case when $V$ is a flat of $M$.  When $A
\subseteq V$, the {\em set} $V-A$ may be considered to be the {\em
matroid} $V/A$ in matroid $M/A$, but considered to be the {\em
matroid} $V-A$ in matroid $M-A$.  We will also use the notions of
internal and external activity as in, \eg, \cite{Bjorner}.

\subsection{A partition}\label{su:partition}

If $\Delta$ is a simplicial complex and $A$ is a set disjoint from the
vertices of $\Delta$, then let $A \ostar \Delta$ denote
$$
A \ostar \Delta := \{A \disun F\colon F \in \Delta\}.
$$
It will soon be important to note that $A \ostar \Delta$ is a
simplicial pair; in fact 
$A \ostar \Delta = (2^A*\Delta, (2^A\setm \{A\})*\Delta)$, 
where $2^A$ denotes the simplicial complex consisting of all subsets
of $A$, and $*$ denotes the usual join, as defined in Section
\ref{se:shifted}.

\begin{lemma}\label{th:ostar}
If $\Delta$ is a simplicial complex and $A$ a finite set disjoint from
the vertices of $\Delta$, then
$$
S_{A \ostar \Delta}(t,q) = t^{\abs{A}}S_\Delta(t,q).
$$
\end{lemma}

\begin{proof}
Under the natural bijection between $\Delta$ and $A \ostar \Delta$,
given by $\phi: F \mapsto A \disun F$, the boundary operators
$\bdmap_\Delta$ and $\bdmap_{A \ostar \Delta}$ are the same.  That is,
$ \bdmap_{A \ostar \Delta}[A \disun F] = [A] \bdmap_\Delta[F]$,
simply by numbering the vertices of $A \ostar \Delta$ so that the
elements of $A$ all come last.  Since the boundary operators are the
same, so are the Laplacians, but the dimension shift in $\phi$ means
$\s_i(\Delta) = \s_{i+\abs{A}}(A \ostar \Delta)$.
The lemma now follows readily.
\end{proof}

If $I$ is independent in $M$ and $p \in \clM{I} - I$,
we will let $\ci(p,I) = \ci_M(p,I) = \ci_{\clM{I}}(p,I)$ be the unique
circuit of $\clM{I}$ contained in $I \disun p$. Dually, if $b \in
I$, we will let $\bo(b,I) = \bo_M(b,I) = \bo_{\clM{I}}(b,I)$ be the
unique bond of $\clM{I}$ contained in $(\clM{I}-I) \disun b$.  
It is easy to see that if $p \not\in \clM{I}$, then $I \in \cI(M/p)$.
Therefore we may safely refer to $\ci_M(p,I)$ for any $I \in \cI(M-p)
- \cI(M/p)$.

\begin{lemma}\label{th:same.ci}
If $I', I \in \cI(M-e) - \cI(M/e)$ and $I' \subseteq I$, then
$\ci_M(e,I') = \ci_M(e,I)$.
\end{lemma}

\begin{proof}
From $\ci_M(e,I') \subseteq I' \disun e \subseteq I \disun e$ it
follows that $\ci_M(e,I')$ is a circuit in $I \disun e$, and thus the
unique circuit in $I \disun e$, \ie, $\ci_M(e,I)$.
\end{proof}

The following lemma is the key step to proving that matroids satisfy
the spectral recursion.

\begin{lemma}\label{th:matroid.partition}
Let $M(E)$ be a matroid, and $e \in E$.  If $e$ is not a loop, then
$$
L_i(M-e,M/e) = 
  \bigoplus_{\substack{C \in \cC(M)\\ e\in C}} 
    L_i((C-e)\ostar \IN(M/C)).
$$
\end{lemma}

\begin{proof}
For any $C \in \cC(M)$ such that $e \in C$, let 
$$
M_C = \{I \in \cI(M-e) - \cI(M/e)\colon \ci_M(e,I)=C \};
$$
we will see shortly that this is a simplicial pair.  By Lemma
\ref{th:same.ci},
$$
\bdmap_{(M-e,M/e)}[I]=\bdmap_C[I]
$$
for any $I \in \cI(M-e) - \cI(M/e)$, where $C=\ci_M(e,I)$.
Thus
removing $M/e$ from $M-e$ partitions
$L_i(M-e,M/e)$ into
$$
L_i(M-e,M/e) = \bigoplus_{\substack{C \in \cC(M)\\ e\in C}} L_i(M_C).
$$
Furthermore, it is easy to see that 
\begin{align*}
M_C 
&= \{I \in \cI(M-e) \colon C-e \subseteq I \}
= (C-e) \ostar \IN((M-e)/(C-e))\\
&= (C-e) \ostar \IN(M/C).
\end{align*}
\end{proof}

\subsection{The Kook-Reiner-Stanton algorithm}\label{su:KRS}
The decomposition in Proposition \ref{th:KRS.algorithm} below was first
discovered by Etienne and Las Vergnas \cite[Theorem
5.1]{EtienneLasVergnas}, but we will rely upon Algorithm \ref{al:KRS},
due to Kook, Reiner, and Stanton \cite[proof of Theorem
1]{KookReinerStanton}, for producing this decomposition.

\begin{proposition}\label{th:KRS.algorithm}
Given a base $B$ of matroid $M$, there is a unique disjoint
decomposition $B=B_1 \disun B_2$ into two (necessarily) independent
sets such that:
\begin{itemize}
\item $B_1$ has internal activity $0$; and
\item $B_2$ has external activity $0$, with respect to the matroid
$M/V$, where $V = \clM{B_1}$.
\end{itemize}
\end{proposition}

\begin{algorithm}\label{al:KRS}
This algorithm produces the decomposition guaranteed by the previous
theorem.  It takes the base $B$ as input, and outputs the pair
$(B_1,B_2)$.
\begin{steplist}
\item Set $B_1=B$, $B_2 = \emptyset$.
\item\label{it:V} Let $V = \clM{B_1}$.
\item Find an internally active element $b$ for $B_1$ as a base of the
flat $V$.
\begin{itemize}
\item If no such element $b$ exists, then stop and output the pair
$(B_1,B_2)$.
\item If such a $b$ exists, then set $B_1 := B_1 - b$, $B_2 := B_2 \disun b$
(we call this step a {\em removal} \!), and return to Step \ref{it:V}. 
\end{itemize}
\end{steplist}
\end{algorithm}

\begin{notation}
If the decomposition of base $B$ in matroid $M$ produced by the above
algorithm is $B = B_1 \disun B_2$, then let $\ppi{}(B)=\ppi{M}(B)
= B_1$.  If $I \in \cI(M)$, then let
$
\pibar{M}(I) = \cl_V(\ppi{V}(I)) = \cl_M(\ppi{V}(I))
$,
where $V=\cl_M(I)$.
If $W$ is any closed set containing $I$ (equivalently, containing
$V=\cl_M(I)$), then $\cl_W(I)=\cl_M(I)=V$, and so
$ \pibar{W}(I) = \cl_V(\ppi{V}(I)) = \pibar{M}(I)$.
In particular,
$\pibar{V}(I)=\pibar{M}(I)$.
\end{notation}

The following lemma, which is little more than a recasting of
\cite[Corollary 18]{KookReinerStanton} in language tailored to our
purposes, reduces computations of the spectrum polynomial to
computations of $\pibar{}$.

\begin{lemma}\label{th:KRS-step}
For any matroid $M(E)$, 
$$
S_M(t,q) 
  = q^{\abs{E}} \sum_{I \in \cI(M)} 
    t^{\rk(\clM{I})} (q^{-1})^{\abs{\pibar{M}(I)}}
  = q^{\abs{E}} \sum_{V \in \cF(M)} 
    t^{\rk(V)} \sum_{I \in \cB(V)}(q^{-1})^{\abs{\pibar{V}(I)}}.
$$
\end{lemma}

Let $\chitcx{\Delta} := \sum (-1)^i f_i(\Delta)$ denote the (reduced)
Euler characteristic of simplicial complex $\Delta$; we also use the
shorthand $\chit{M}=\chit{\IN(M)}$. If $V \subseteq W$ are flats of
matroid $M$, let $\mu(W,V)=\mu_M(W,V)$ denote the M\"obius
function of the sublattice $[W,V]$ in the lattice of flats of $M$.
The proof of \cite[equation (2.2)]{KookReinerStanton} shows that
\begin{equation}\label{eq:KRS-step}
\sum_{B \in \cB(M)} x^{\abs{\pibar{M}(B)}}
 = \sum_{V \in \cF(M)} \abs{\chit{V}} \abs{\mu(V,M)} x^{\abs{V}}.
\end{equation}
We use the same techniques to
do something similar.

\begin{lemma}\label{th:e-step}
For any matroid $M(E)$, and any $e \in E$,
$$
\sum_{\substack{B \in \cB(M)\\e \in \pibar{M}(B)}} 
  x^{\abs{\pibar{M}(B)}}
 = \sum_{\substack{V \in \cF(M) \\ e \in V}} 
  \abs{\chit{V}} \abs{\mu(V,M)} x^{\abs{V}}.
$$
In particular, this sum is independent of the linear order on $E$.
\end{lemma}

\begin{proof}
By Algorithm \ref{al:KRS} (see also its proof in \cite{KookReinerStanton}),
there is a bijection between:
\begin{itemize}
\item 
the set $\cV$ of triples $(V,B_1,B_2)$ where $V$ is a flat of
$M$, $B_1$ is a base of internal activity 0 for $V$ (in particular,
$V = \clM{B_1}$), and $B_2$ is a base of
external activity 0 for $M/V$; and
\item
the set $\cB$ of bases $B$ of $M$.
\end{itemize}
Furthermore, $B=B_1 \disun B_2$ and $\ppi{M}(B)=B_1$.
Thus
$$
\sum_{\substack{B \in \cB(M)\\e \in \pibar{M}(B)}} 
  x^{\abs{\pibar{M}(B)}}
 = \sum_{\substack{(V,B_1,B_2) \in \cV\\e \in V}} x^{\abs{V}}.
$$

We must then determine how many triples $(V,B_1,B_2)$ there are in
$\cV$ for a fixed flat $V$.  Mimicking an argument from the proof of
\cite[Theorem 1]{KookReinerStanton}, we recall from \cite[Theorem
7.8.4]{Bjorner} that there are $\abs{\chit{V}}$ bases of internal
activity 0 for $V$, and from \cite[Proposition 7.4.7]{Bjorner} that
there are $\abs{\mu(V,M)}$ bases of external activity 0 for $M/V$.  So
for every $V$, there are $\abs{\chit{V}}$ choices for $B_1$, and,
independently, $\abs{\mu(V,M)}$ choices for $B_2$.  Thus,
$$
\sum_{\substack{(V,B_1,B_2) \in \cV\\e \in V}} x^{\abs{V}}
 = \sum_{\substack{V \in \cF(M)\\e \in V}} 
  \abs{\chit{V}} \abs{\mu(V,M)} x^{\abs{V}},
$$
completing the proof.
\end{proof}

We now see how Algorithm
\ref{al:KRS} works on $M-e$ (Lemma \ref{th:mat.del}) and $M/e$ (Lemma
\ref{th:mat.lk}), and on $M/C$ when $C$ is a
circuit containing $e$ (Lemma \ref{th:KRS.circuit}). We first need
three technical lemmas whose easy proofs are omitted.
We abuse set difference notation slightly
to let $A \setm x$ denote $\{a \in A \colon a \neq x\}$, when $A$ is a
set that may or may not contain element $x$.

\begin{lemma}\label{th:A3}
Let $I$ be an independent set in matroid $M$, let $e$ be last in the
linear order, and assume that $e \not\in I$ and that $e$ is not an
isthmus of $M$.  Then $b$ is internally active in $I$ (with respect to
$M$) iff $b$ is internally active in $I$ (with respect to $M-e$).
\end{lemma}

\begin{lemma}\label{th:A1}
Let $I$ be an independent set in matroid $M$, and let $e,b \in I$.
Then $b$ is internally active in $I$ (with respect to $M$) iff $b$ is
internally active in $I-e$ (with respect to $M/e$).
\end{lemma}

\begin{lemma}\label{th:A2}
Let $I$ be an independent set in matroid $M$, and let $i$ be an
isthmus in $\clM{I}$.  Then $b \neq i$ is internally active in $I$
(with respect to $M$) iff $b$ is internally active in $I-i$ (with
respect to $M$).
\end{lemma}

\begin{lemma}\label{th:mat.del}
Let $B$ be a base of $M-e$, so $B$ is also a base of $M$ and $e
\not\in B$.  Also assume $e$ is last in the linear order. Then
$\ppi{M-e}(B) = \ppi{M}(B)$.
\end{lemma}

\begin{proof}
Use Algorithm \ref{al:KRS} to compute $\ppi{M}(B)$.  By Lemma
\ref{th:A3}, every step of the algorithm can be copied in $M-e$; that
is, when element $b$ is removed from $B_1$ in $M$, we can remove $b$
from $B_1$ in $M-e$.  And also by Lemma \ref{th:A3}, when there are no
more elements to remove from $B_1$ in $M$, then there are also no more
elements to remove from $B_1$ in $M-e$.
\end{proof}

\begin{corollary}\label{th:mat.del.cor}
Let $I$ be an independent set of $M-e$, so $I$ is also independent in
$M$ and $e \not\in I$.  Also assume $e$ is last in the linear order. Then
$$
\pibar{M-e}(I) = \pibar{M}(I) \setm e.
$$
\end{corollary}

\begin{lemma}\label{th:mat.lk}
Let $B$ be a base of $M$ such that $e \in B$, so $B-e$ is a base of
$M/e$.  Also assume $e$ is last in the linear order. Then
$$\ppi{M/e}(B-e) = 
\ppi{M}(B) \setm e
$$
\end{lemma}

\begin{proof}
Again use Algorithm \ref{al:KRS} to compute $\ppi{M}(B)$, except do
not remove $e$ unless it is the only element that can be removed.  As
in Lemma \ref{th:mat.del}, every step can be copied in $M/e$, this
time by Lemma \ref{th:A1}, as long as we are not removing $e$, and
have not yet removed $e$.  Also by Lemma \ref{th:A1}, if we never
remove $e$, then when there are no more elements to remove in $M$,
there are no more elements to remove in $M/e$.  Thus, if $e$ is never
removed (\ie, if $e \in \ppi{M}(B))$, then $\ppi{M/e}(B-e) =
\ppi{M}(B)-e$.

If $e$ is eventually removed in $M$, it must be when $e$ is an
isthmus, since $e$ is ordered last (so it can be the minimal
element of $\bo(e,I)$ only if it is the only element -- \ie, if it is an
isthmus).  Since we put off removing $e$ until there were no other
possible removals, Lemma \ref{th:A2} guarantees that there are no new
removals possible after $e$ is removed.  Since the removals were
identical in $M$ and $M/e$ until $e$ was removed, 
$\ppi{M/e}(B-e) = \ppi{M}(B)$.
\end{proof}

\begin{corollary}\label{th:mat.lk.cor}
Let $I$ be an independent set of $M$ such that $e \in I$, so $I-e$ is
independent in $M/e$.  Also assume $e$ is last in the linear order. Then
$$
\pibar{M/e}(I-e) = \pibar{M}(I) \setm e.
$$
\end{corollary}

\begin{proof}
Let $V = \cl_M(I)$.  Then 
$\cl_{M/e}(I-e)= V-e$ 
as {\em sets}, so $\cl_{M/e}(I-e) = V/e$ as {\em matroids}.  Thus, by
the definition of $\pibar{}$, we have
$$
\pibar{M/e}(I-e) =\cl_{M/e}(\ppi{V/e}(I-e)).
$$

If $e \in \ppi{V}(I)$, then simply
$
\cl_{M/e}(\ppi{V/e}(I-e))
 =\cl_{M/e}(\ppi{V}(I)-e)
  =\cl_{M}(\ppi{V}(I))-e 
 =\pibar{M}(I) \setm e;
$
the first equality is by Lemma \ref{th:mat.lk}, the second equality is
a routine exercise using $e \in \ppi{V}(I)$, and the last equality is
from the definition of $\pibar{}$.

If $e \not\in \ppi{V}(I)$, then the proof of Lemma \ref{th:mat.lk}
shows that $e$ is an isthmus in $\cl_V(\ppi{V}(I)\cup e)$.  Then, since
$\cl(A \disun i) = (\cl A) \disun i$ for any $A$ and any isthmus $i
\not\in A$,
\begin{equation}\label{eq:mat.lk.cl}
\cl_M(\ppi{V}(I)\cup e) = \cl_V(\ppi{V}(I)\cup e) = \cl_V(\ppi{V}(I))\cup e = 
  \pibar{M}(I) \cup e.
\end{equation}
Now, also in this case,
\begin{align*}
\cl_{M/e}(\ppi{V/e}(I-e))
 &=\cl_{M/e}(\ppi{V}(I))
  =\cl_{M}(\ppi{V}(I) \cup e)-e
  =(\pibar{M}(I) \cup e)-e\\
 &=\pibar{M}(I) \setm e;
\end{align*}
the first equality is by Lemma \ref{th:mat.lk}, the second equality is from
the definition of $\cl_{M/e}$, and the third equality is equation
\eqref{eq:mat.lk.cl}.
\end{proof}

\begin{lemma}\label{th:KRS.circuit}
Let $B$ be a base of matroid $M(E)$, let $e$ be first in the linear
order on $E$, and assume that $e \not\in B$ and $e$ is not a loop.
Let $C=\ci(e,B)$, so $B-(C-e)$ is a base of $M/C$.  Then
$$
\ppi{M/C}(B-(C-e)) = \ppi{M}(B) - (C-e).
$$
\end{lemma}

\begin{proof}
It is an easy exercise to check
that $\bo_{M/C}(b,B-(C-e))=\bo_M(b,B)$ for any $b \in
B-(C-e)$.  
It then follows 
that $b$ is internally active in $B$
(with respect to $M$) iff $b$ is minimal in
$\bo_M(b,B)=\bo_M(b,B-(C-e))$ iff $b$ is internally active in
$B-(C-e)$ (with respect to $M/C$).

Now, as 
in Lemmas \ref{th:mat.del} and \ref{th:mat.lk},
use Algorithm \ref{al:KRS} to compute
$\ppi{M/C}(B-(C-e))$.  Once again, every step can be copied in $M$,
computing $\ppi{M}(B)$.  Furthermore, when there are no more elements
in $B-(C-e)$ to remove in computing $\ppi{M/C}(B-(C-e))$, the only
elements of $B$ that could possibly be removed in computing
$\ppi{M}(B)$ must be in $C-e$.  We now show that any $c \in C-e$ is
not internally active, and thus that the removals in $M$ and $M/C$ are
identical, which will complete the proof.

It is easy to see that $C = \ci_{\clM{B_1}}(e,B_1)$, where $B_1$ is
what remains of $B$ after performing all the removals in $M$
corresponding to the removals in $M/C$.  Thus $c \in C-e \subseteq
\ci_{\clM{B_1}}(e,B_1)$ implies, by 
\eg, \cite[Lemma 7.3.1]{Bjorner},
that
$e \in \bo_{\clM{B_1}}(c,B_1)$.  Since $e$ is first in the linear
order, $c$ is, as desired, not internally active.
\end{proof}

\subsection{The spectral recursion for matroids}\label{su:mat.finish}

We now prove that matroids satisfy the spectral recursion (Theorem
\ref{th:matroids}), by comparing $qtS_{M/e} + qS_{M-e} - S_M$ and
$S_{(M-e,M/e)}$.
In each
case, we get two expressions, one in terms of $\chitempty$ and $\mu$, the
other in terms of $\pibar{}$.  The expressions in terms of $\chitempty$ and $\mu$
lead to a quick proof, by reducing a key piece of the equation to the
$q=0$ case for a flat.  The expressions in terms of $\pibar{}$ suggest a
more bijective proof, which is not hard to prove either.  Both proofs
are given.

\begin{lemma}\label{th:p9}
If $M(E)$ is a matroid, and $e \in E$ is neither an isthmus nor a loop, then
\begin{multline*}
q S_{M-e}(t,q) + qt S_{M/e}(t,q) - S_M(t,q)\\
 \begin{aligned}
  &= (q-1) q^{\abs{E}}  \sum_{V \in \cF(M)} t^{\rk_M(V)}
     \sum_{\substack{I \in \cB(V)\\e \in \pibar{V}(I)}}
     (q^{-1})^{\abs{\pibar{V}(I)}}\\
  &= (q-1) \sum_{V \in \cF(M)} t^{\rk_M(V)}
     \sum_{\substack{W \in \cF(V)\\e \in W}}
     \abs{\chit{W}} \abs{\mu(W,V)} q^{\abs{E}-\abs{W}}.
 \end{aligned}
\end{multline*}
\end{lemma}

\begin{proof}
We compute each of $S_{M-e}$ and $S_{M/e}$ using Lemma \ref{th:KRS-step}.
First,
\begin{align}
S_{M-e}(t,q) 
 &= q^{\abs{E-e}} \sum_{I \in \cI(M-e)}
      t^{\rk_{M-e}(\cl_{M-e}(I))} (q^{-1})^{\abs{\pibar{M-e}(I)}} \nonumber \\
 &= q^{\abs{E-e}} \sum_{\substack{I \in \cI(M)\\e \not\in I}}
      t^{\rk_M(\cl_M(I))}(q^{-1})^{\abs{\pibar{M}(I) \setm e}},\label{eq:S.M-e}
\end{align}
since: $I \in \cI(M-e)$ iff $I \in \cI(M)$ and $e \not\in I$;
$\abs{\pibar{M-e}(I)} = \abs{\pibar{M}(I) \setm e}$, by Corollary
\ref{th:mat.del.cor}; and $\rk_{M-e}(\cl_{M-e}(I)) =
\rk_{M}(\cl_{M}(I)\setm e)$ is an easy matroid exercise.
Similarly,
\begin{align}
S_{M/e}(t,q)
 &= q^{\abs{E-e}} \sum_{I' \in \cI(M/e)}
      t^{\rk_{M/e}(\cl_{M/e}(I'))} (q^{-1})^{\abs{\pibar{M/e}(I')}} \nonumber \\
 &= q^{\abs{E-e}} \sum_{\substack{I \in \cI(M)\\e \in I}}
      t^{\rk_M(\cl_M(I))-1}(q^{-1})^{\abs{\pibar{M}(I) \setm e}}, \label{eq:S.M/e}
\end{align}
where
$I = I' \disun e$ for $I' \in \cI(M/e)$, since:
$\abs{\pibar{M/e}(I')} =
\abs{\pibar{M/e}(I-e)} = \abs{\pibar{M}(I) \setm e}$, by Corollary \ref{th:mat.lk.cor}; and
$
\rk_{M/e}(\cl_{M/e}(I')) 
  = \rk_M(\cl_M(I))-1
$
is a routine exercise,
using $e \in \cl_M(I)$.

Combining equations \eqref{eq:S.M-e} and \eqref{eq:S.M/e}, and then
sorting independent sets by their closures, we get
\begin{align}
q S_{M-e}(t,q) + qt S_{M/e}(t,q) 
  &=q^{\abs{E}} \sum_{I \in \cI(M)}
     t^{\rk_M(\cl_M(I))}(q^{-1})^{\abs{\pibar{M}(I) \setm e}} \nonumber \\
  &= q^{\abs{E}} \sum_{V \in \cF(M)}t^{\rk_M(V)}
     \sum_{I \in \cB(V)} (q^{-1})^{\abs{\pibar{V}(I) \setm e}}.\label{eq:S.M-e.M/e}
\end{align}
Furthermore
\begin{align}
\sum_{I \in \cB(V)} (q^{-1})^{\abs{\pibar{V}(I) \setm e}}
  &=  \sum_{\substack{I \in \cB(V)\\e \not\in \pibar{V}(I)}} 
	(q^{-1})^{\abs{\pibar{V}(I) \setm e}}
    +\sum_{\substack{I \in \cB(V)\\e \in \pibar{V}(I)}} 
	(q^{-1})^{\abs{\pibar{V}(I) \setm e}}\nonumber \\
  &= \sum_{I \in \cB(V)} 
	(q^{-1})^{\abs{\pibar{V}(I)}}
     +(q-1) \sum_{\substack{I \in \cB(V)\\e \in \pibar{V}(I)}}
	 (q^{-1})^{\abs{\pibar{V}(I)}}; \label{eq:q.piece}
\end{align}
plugging into equation \eqref{eq:q.piece} into equation \eqref{eq:S.M-e.M/e} readily leads to the first equation
of the lemma.
The second equation then follows directly from Lemma \ref{th:e-step}.
\end{proof}

\begin{lemma}\label{th:p11}
If $M(E)$ is a matroid, and $e \in E$ is neither an isthmus nor a
loop, then
\begin{multline*}
S_{(M-e,M/e)}(t,q)
  = q^{\abs{E}}\sum_{V \in \cF(M)} t^{\rk_M(V)}
      \sum_{\substack{C \in \cC(V)\\e \in C}} q^{-\abs{C}}
        \sum_{I \in \cB(V/C)}(q^{-1})^{\abs{\pibar{V/C}(I)}}\\
  = \sum_{V \in \cF(M)} t^{\rk_M(V)}
	\sum_{\substack{W \in \cF(V)\\e \in W}} 
	\sum_{\substack{C \in \cC(W)\\e \in W}} 
              \abs{\chit{W/C}}
	\abs{\mu(W,V)} q^{\abs{E}-\abs{W}}.
\end{multline*}
\end{lemma}

\begin{proof}
By Lemmas \ref{th:KRS-step}, \ref{th:ostar}, and \ref{th:matroid.partition},
\begin{align*}
S_{(M-e,M/e)}(t,q)
 &=\sum_{\substack{C \in \cC(M)\\e \in C}} t^{\rk_M(C)}S_{M/C}(t,q)\\
 &=\sum_{\substack{C \in \cC(M)\\e \in C}} t^{\rk_M(C)} 
    q^{\abs{E-C}} \sum_{W \in \cF(M/C)} t^{\rk_{M/C}(W)}
      \sum_{I \in \cB(W)} (q^{-1})^{\abs{\pibar{W}(I)}}.
\end{align*}
Now, the flats of $M/C$ are $V-C$ as {\em sets}, and thus $V/C$ as
{\em matroids}, for all flats $V$ of $M$ containing $C$.  Therefore, 
\begin{multline*}
S_{(M-e,M/e)}(t,q)\\
\begin{aligned}
 &= q^{\abs{E}} \sum_{\substack{C \in \cC(M)\\e \in C}} t^{\rk_M(C)} 
    q^{-\abs{C}} \sum_{\substack{V \in \cF(M)\\C \subseteq V }} 
        t^{\rk_{M}(V)-\rk_M(C)}
      \sum_{I \in \cB(V/C)} (q^{-1})^{\abs{\pibar{V/C}(I)}}\\
 &= q^{\abs{E}} \sum_{V \in \cF(M)} t^{\rk_M(V)}
     \sum_{\substack{C \in \cC(M)\\e \in C\subseteq V}} q^{-\abs{C}}
     \sum_{I \in \cB(V/C)} (q^{-1})^{\abs{\pibar{V/C}(I)}},
\end{aligned}
\end{multline*}
which is the first equation of the lemma, once we note that $C \in \cC(V)$ iff $C \in \cC(M)$ and $C \subseteq V$.  

The second
equation of the lemma then follows from
\begin{multline*}
\sum_{\substack{C \in \cC(V)\\e \in C}} q^{-\abs{C}}
      \sum_{I \in \cB(V/C)} (q^{-1})^{\abs{\pibar{V/C}(I)}}\\
 \begin{aligned}
  &=\sum_{\substack{C \in \cC(V)\\e \in C}} q^{-\abs{C}}
	\sum_{W/C \in \cF(V/C)}
	  \abs{\chit{W/C}} \abs{\mu_{V/C}(W/C,V/C)} (q^{-1})^{\abs{W/C}}\\
  &=\sum_{\substack{C \in \cC(V)\\e \in C}} 
	\sum_{\substack{W \in \cF(V)\\C \subseteq W}}
	  \abs{\chit{W/C}} \abs{\mu(W,V)} (q^{-1})^{\abs{W}}\\
  &=\sum_{\substack{W \in \cF(V)\\e \in W}}
	\sum_{\substack{C \in \cC(V)\\e \in C \subseteq W}}
	  \abs{\chit{W/C}} 
        \abs{\mu(W,V)} (q^{-1})^{\abs{W}}.
 \end{aligned}
\end{multline*}
The first equation above is from equation \eqref{eq:KRS-step}; we are
also using the same characterization of flats of a contraction as in
the previous paragraph.  The second equation is since the interval
$[W/C,V/C]$ in the lattice of flats of $V/C$ is isomorphic to the
interval $[W,V]$ in the lattice of flats of $V$, again by that same
characterization of flats in a contraction.  It only remains to again
note that $C \in \cC(W)$ iff $C \in \cC(M)$ and $C \subseteq W$.
\end{proof}

\begin{theorem}\label{th:matroids}
If $M$ is a matroid, then its independence complex $\IN(M)$ satisfies
the spectral recursion, equation \eqref{eq:big}.
\end{theorem}

\begin{proof}
By Proposition \ref{th:loop.ok}, we may assume $e$ is not a loop.  By Lemma
\ref{th:v.ok} and Theorem \ref{th:join.ok.e} below (which does not
depend on anything in this section), we may assume $e$ is not an
isthmus. As discussed at the beginning of the subsection, there are
now two ways to finish off the proof, one using the $q=0$ case, the
other using a bijection.

{\em $q=0$ proof.}  
By Theorem \ref{th:homology}, we know that the spectral recursion holds, for
any matroid, with $q=0$.  By Lemmas \ref{th:p9} and \ref{th:p11}, this
means
\begin{equation}\label{eq:magic.chi}
\abs{\chit{M}} = \sum_{\substack{C \in \cC(M)\\e \in C}}
	\abs{\chit{M/C}}
\end{equation}
for any matroid $M$, since only terms with $W=E$ survive when $q=0$.
(Equation \eqref{eq:magic.chi} is also, as noted by Kook
\cite{Kook:circuit}, dual to Crapo's complementation theorem (\eg,
\cite[Theorem 4.33]{Aigner}) applied to the dual matroid of $M$.)
Thus, simply by plugging in the flat $W$, as a matroid, for the
matroid $M$ in equation \eqref{eq:magic.chi},
$$
\abs{\chit{W}} = \sum_{\substack{C \in \cC(W)\\e \in C}}
	\abs{\chit{W/C}}
$$
whenever $W$ is a flat of $M$ containing $e$.
By Lemmas \ref{th:p9} and \ref{th:p11} again, we are done.

{\em Bijective proof.}
By Lemmas \ref{th:p9} and \ref{th:p11}, it suffices to show
\begin{equation}\label{eq:more.bijective}
\sum_{\substack{I \in \cB(M)\\e \in \pibar{W}(I)}}(q^{-1})^{\abs{\pibar{M}(I)}}
= \sum_{\substack{C \in \cC(M)\\e \in C}}
    \sum_{I \in \cB(M/C)} (q^{-1})^{\abs{\pibar{M/C}(I)}+\abs{C}}.
\end{equation}
Further, Lemma \ref{th:e-step} shows that the sum on the left-hand side
of equation \eqref{eq:more.bijective} is independent of the ordering
of the ground set.  Similarly, Lemma \ref{th:p11} itself shows the
same thing for the sum on the right-hand side.  So we now assume, for
the remainder of this proof, that $e$ is ordered first in the linear
order on $E$.

Equation \eqref{eq:more.bijective} would follow naturally from a bijection
$$
\phi\colon \{B \in \cB(M)\colon e \in \pibar{M}(B)\}
 \rightarrow \{(C,I)\colon C \in \cC(M), I \in \cB(M/C), e \in C\}
$$
such that 
\begin{equation}\label{eq:pi.bijection}
\pibar{M/C}(I) \disun C = \pibar{M}(B),
\end{equation}
where $\phi(B)=(C,I)$.  Such a bijection is given by, as we now show,
$C = \ci(e,B)$ and $I=B-(C-e)$ in one direction, and $B=I \disun C -
e$ in the other.

First note that, since $e$ is ordered first, if $e \in B$ then $e$ is
internally active in $B$, and so $e \not\in \ppi{M}(B)$.  It is then
easy to see in this case that $e \not\in \pibar{M}(B)$.  We may
therefore safely assume $e \not\in B$, and so $C = \ci(e,B)$ is
well-defined.  It then follows that $\phi$ is well-defined.

It is easy to see that $\phi$ is injective.  Showing that $\phi$ is
surjective reduces to verifying that $e \in \pibar{M}(B)$ when $B=I
\disun C - e$; by Lemma \ref{th:KRS.circuit}, $C-e \subseteq
\ppi{M}(B)$, so $e \in C = \clM{C-e} \subseteq \pibar{M}(B)$.

Finally, to verify equation \eqref{eq:pi.bijection}, by Lemma
\ref{th:KRS.circuit} and the definition of closure in a matroid
contraction,
$
\pibar{M/C}(B-(C-e))
 = \cl_M(\ppi{M}(B) \cup e) - C.
$
Also, $\cl_M(\ppi{M}(B)\cup e) = \cl_M(\ppi{M}(B)) = \pibar{M}(B)$,
since $e \in \pibar{M}(B)$, which completes the proof of equation
\eqref{eq:pi.bijection}.
\end{proof}

\begin{remark}\label{rm:truly.recursive}
The spectral recursion does not provide a truly recursive way to
compute $S_M$, due to the presence of $S_{(M-e,M/e)}$, since the
recursion only applies to a single matroid, and not a matroid pair
like $(M-e,M/e)$.  We can however, combine it with Lemmas
\ref{th:ostar} and \ref{th:matroid.partition} for a recursion that is
truly recursive, albeit with more terms than the spectral recursion:
$$
S_M(t,q) = qS_{M-e}(t,q) + qtS_{M/e}(t,q)
      + (1-q)\sum_{\substack{C \in \cC(M)\\e \in C}} t^{\rk_M(C)}S_{M/C}(t,q).
$$
I am grateful to E. Babson for this observation.
\end{remark}

\section{Shifted complexes}\label{se:shifted}
We postpone until Subsection \ref{su:shift.finish} the actual
definition of shifted complexes, but we will see there that a shifted
complex is a skeleton of a cone of a smaller shifted complex (Lemmas
\ref{th:inductive-shifted-definition} \ref{th:pure.shifted}.  To prove
that shifted complexes satisfy the spectral recursion, equation
\eqref{eq:big}, then, it suffices to show that taking skeleta and
taking cones each preserve the property of satisfying the spectral
recursion -- which are interesting results in their own right.

We will prove in Subsection \ref{su:joins} that the property of
satisfying the spectral recursion is preserved by taking joins
(Corollary \ref{th:join.ok.all}), and thus by taking cones (\cf\
Proposition \ref{th:v.ok}).  The key step is that a simple formula
\cite[Theorem 4.10]{DuvalReiner} for the eigenvalues of the join
generalizes straightforwardly from single simplicial complexes to
simplicial pairs (Corollaries \ref{spectra-join} and \ref{th:S.join}).

Proving that taking skeleta preserves the property of satisfying the
spectral recursion is harder, and is the focus of Subsections
\ref{su:fine.laplace}--\ref{su:spectral.skeleta}.  The key facts about
Laplacians, established in Subsections \ref{su:fine.laplace} and
\ref{su:skeleta}, respectively, are that the non-zero eigenvalues come
in pairs in consecutive dimensions (Lemma \ref{th:precursor}), and
that taking $(d-1)$-skeleta preserves non-zero eigenvalues of the
finer Laplacians in dimension $d-1$ and below (Lemma
\ref{th:little.s.d-1}).

The only eigenvalues in dimension $d-1$ and below that are changed by
taking $(d-1)$-skeleta, then, are some $(d-1)$-dimensional eigenvalues
that become $0$ when their counterparts (in the sense of Lemma
\ref{th:precursor}) in dimension $d$ are removed.  It is auspicious
that these replaced $(d-1)$-dimensional eigenvalues must line up
properly in the spectral recursion (since their counterparts in dimension
$d$, the only non-zero eigenvalues in that dimension, do as well) and
that the $0$'s that replace them also line up properly (since the spectral
recursion is true with $q=0$ for both the original complex and its
skeleton, by Theorem \ref{th:homology}).  But it turns out that we are
better off with $f$-vectors ($q=1$, also a good case by Theorem
\ref{th:homology}) than with homology ($q=0$), in part because the
change in $f$-vectors resulting from taking skeleta is much easier to
describe than the change in homology.

In Subsection \ref{su:spectral.skeleta}, we will see that the
difference between the spectrum polynomials of the skeleton and the
original complex can be described largely in terms of the $f$-vector
(Lemma \ref{th:S.skeleton}), allowing us to describe the difference in
the spectral recursion between the skeleton and the original complex
in a particularly useful form (Lemma \ref{th:scriptS.skel}).  From
there, simple generating function manipulations lead to Theorem
\ref{th:skeleta.preserve}, which states that a $d$-dimensional
simplicial complex satisfies the spectral recursion with respect to a
vertex if and only if its $(d-1)$-skeleton and pure $d$-skeleton (the
complex generated by its facets) do as well.

\subsection{Joins and cones}\label{su:joins}
Define the {\em join} $(\Delta,\Delta') *
(\Gamma,\Gamma')$ of two simplicial pairs on disjoint vertex sets to
be
$$
(\Delta,\Delta') * (\Gamma,\Gamma') := 
\{F \disun G \colon F \in \Delta \setm \Delta', G \in \Gamma \setm \Gamma'\}
$$
(here, $\disun$ denotes disjoint union),
which equals the simplicial pair
\begin{equation}\label{eq:nice.join}
(\Delta*\Gamma, (\Delta'*\Gamma) \cup (\Delta*\Gamma')).
\end{equation}
When $\Delta'=\Gamma'=\emptyset$, this reduces to the usual join
$\Delta*\Gamma$.  When, further, $\Delta$ is a single vertex, say $v$,
the join is written as $v*\Gamma$, the {\em cone} over $\Gamma$ with
{\em apex} $v$.

The proofs of the following two results on simplicial pairs are
identical (modulo some indexing changes) to those of the analogous
statements for single simplicial complexes \cite[Section
4]{DuvalReiner}.

\begin{proposition}
For any two simplicial pairs $(\Delta, \Delta')$ and $(\Gamma,
\Gamma')$ and every $k$, the map defined $\reals$-linearly by
$
[F] \otimes [G] \mapsto [F \disun G]
$
identifies the vector spaces 
$$
\bigoplus_{i+j=k} 
C_{i-1}((\Delta,\Delta');\reals) \otimes 
   C_{j-1}((\Gamma,\Gamma');\reals)
\cong C_{k-1}((\Delta,\Delta') * (\Gamma,\Gamma') ; \reals)
$$
and has the following property with respect to the Laplacians $L$ of
the appropriate dimensions in 
$(\Delta,\Delta'), (\Gamma,\Gamma')$, and 
$(\Delta,\Delta') * (\Gamma,\Gamma')$:
\begin{equation}
\label{join-operator-equation}
L((\Delta,\Delta') * (\Gamma,\Gamma')) 
    = L(\Delta,\Delta') \otimes \id 
      + \id \otimes L(\Gamma,\Gamma').
\end{equation}
\end{proposition}

\begin{corollary}\label{spectra-join}
If $(\Delta, \Delta')$ and $(\Gamma, \Gamma')$ are two simplicial
pairs, 
then
$$
\s_{k-1}((\Delta, \Delta') *(\Gamma, \Gamma') ) = 
  \bigcup_{\condns{i+j=k}{\lambda 
		\in \s_{i-1}(\Delta, \Delta'),\ 
            \mu \in \s_{j-1}(\Gamma, \Gamma')}} 
    \lambda+\mu.
$$
\end{corollary}

It is then an easy exercise in generating functions to verify the
following corollary.

\begin{corollary}\label{th:S.join}
If $(\Delta, \Delta')$ and $(\Gamma, \Gamma')$ are two simplicial
pairs, 
then
$$
S_{(\Delta, \Delta') * (\Gamma, \Gamma')} = 
 S_{(\Delta, \Delta')} S_{(\Gamma, \Gamma')}.
$$
\end{corollary}

\begin{theorem}\label{th:join.ok.e}
If $\Delta$ satisfies the spectral recursion with respect to $e$, and 
$\Gamma$ is any simplicial complex whose vertex set is disjoint from 
the vertex set of $\Delta$, then the join $\Delta * \Gamma$ satisfies the spectral 
recursion with respect to $e$.    
\end{theorem}
\begin{proof}
By Corollary \ref{th:S.join} twice, and our hypothesis,
\begin{align*}
S_{\Delta * \Gamma} = S_{\Delta} S_{\Gamma}
  &= (qS_{\Delta-e} + qtS_{\Delta/e} + (1-q)S_{(\Delta-e,\Delta/e)})S_{\Gamma}\\
  &= qS_{(\Delta-e)*\Gamma} + qtS_{(\Delta/e)*\Gamma} + (1-q)S_{(\Delta-e,\Delta/e)*\Gamma}.
\end{align*}
This last expression is exactly what we need, since it is easy to
verify that join commutes with deletion and contraction, \ie, 
$(\Delta-e)*\Gamma = (\Delta*\Gamma)-e$ and
$(\Delta/e)*\Gamma/e = (\Delta*\Gamma)/e$,
and also since equation \eqref{eq:nice.join} with $\Delta'=\emptyset$ then yields
$$
(\Delta-e,\Delta/e)*\Gamma 
 = ((\Delta-e)*\Gamma,(\Delta/e)*\Gamma) 
 = ((\Delta*\Gamma)-e,(\Delta*\Gamma)/e).
$$
\end{proof}

\begin{corollary}\label{th:join.ok.all}
If $\Delta$ and $\Gamma$ each satisfy the spectral recursion, then so
does their join $\Delta*\Gamma$.
\end{corollary}

\subsection{Finer Laplacians}\label{su:fine.laplace}
Recall from Section \ref{se:laplace} that
$L'_i = L'_i(\Delta,\Delta') := \bdmap_{i+1} \bdmap_{i+1}^*$ and
$L''_i = L''_i(\Delta,\Delta') := \bdmap_{i}^* \bdmap_{i},$
so that
$L_i = L'_i + L''_i$.
Define $\s'_i(\Delta,\Delta')$ and $\s''_i(\Delta,\Delta')$ to be the
multiset of eigenvalues of $L'_i(\Delta,\Delta')$ and
$L''_i(\Delta,\Delta')$, respectively, arranged in weakly decreasing order.

Following \cite{DuvalReiner}, let the equivalence relation $\blambda
\circeq \bmu$ on multisets $\blambda$ and $\bmu$ denote that
$\blambda$ and $\bmu$ agree in the multiplicities of all of their {\em
non-zero parts}, \ie, that they coincide except for possibly their
number of zeroes.  Also let $\blambda \cup \bmu$ denote the
$\circeq$-equivalence class whose non-zero parts are the multiset
union of the non-zero parts of $\blambda$ and $\bmu$.

\begin{proposition}\label{th:s.s''}
If $(\Delta,\Delta')$ is a simplicial pair, then
$$
\s_i(\Delta,\Delta') \circeq \s''_i(\Delta,\Delta') \cup \s''_{i+1}(\Delta,\Delta').
$$
\end{proposition}
\begin{proof}
The proof is identical to the single simplicial complex
($\Delta'=\emptyset$) case in \cite[Equation (3.6)]{DuvalReiner}, and
depends only upon $\bdmap^2=0$ and routine eigenvalue calculations
involving adjoints.
\end{proof}

If $(\Delta,\Delta')$ is a simplicial pair, let
\begin{gather*}
S''_{(\Delta,\Delta'),i}(q) := 
  \sum_{\substack{\lambda \in \s''_i(\Delta,\Delta') \\ \lambda \neq 0 }}
    q^\lambda, \quad \text{and}\\
S''_{(\Delta,\Delta')}(t,q) := 
  \sum_{i} S''_{(\Delta,\Delta'),i-1}(q)t^i.
\end{gather*}
Zero eigenvalues are omitted from these definitions of $S''$ in order
to more naturally encode Proposition \ref{th:s.s''} into the language
of generating functions, in Lemma \ref{th:precursor}, below.
Also let 
$$
B_{(\Delta,\Delta')}(t) 
 := \sum_{i} \tilde{\beta}_{i-1}(\Delta,\Delta')t^i
  = \sum_{i} m_0(L_{i-1}(\Delta,\Delta')) t^i
  = S_{(\Delta,\Delta')}(t,0).
$$
These three definitions of $B$ are equivalent by Proposition
\ref{th:beta.m0}.

From now on, when there is no confusion about the variables $t$ and
$q$, we will often omit them for clarity.

\begin{lemma}\label{th:precursor}
If $(\Delta,\Delta')$ is a simplicial pair, then
$$
S_{(\Delta,\Delta')} = (1+t^{-1})S''_{(\Delta,\Delta')} + B_{(\Delta,\Delta')}.
$$
\end{lemma}
\begin{proof}
Combine Propositions \ref{th:s.s''} and \ref{th:beta.m0}.
\end{proof}

\begin{corollary}\label{th:t-1}
If $\Delta$ is any simplicial complex, and $e$ is any vertex of
$\Delta$, then the spectral recursion holds when $t=-1$.
\end{corollary}
\begin{proof}
By Lemma \ref{th:precursor}, for any simplicial pair
$(\Delta,\Delta')$,
$$
S_{(\Delta,\Delta')}(-1,q) = B_{(\Delta,\Delta')}(-1,q)
= \sum_i (-1)^i \tilde{\beta}_i(\Delta,\Delta') =\chi(\Delta,\Delta'),
$$
where $\chi(\Delta,\Delta')$ denotes the Euler characteristic of the
simplicial pair $(\Delta,\Delta')$ (see \eg, \cite{Munkres}).  The
identity $\chi(\Delta,\Delta')=\chi(\Delta)-\chi(\Delta')$, which
holds as long as $\Delta' \subseteq \Delta$, immediately reduces the
$t=-1$ instance of the spectral recursion to
$\chi(\Delta)=\chi(\Delta-e)-\chi(\Delta/e)$.  This, in turn, follows
from $\chi(\Delta)=\sum_i(-1)^i f_i(\Delta)$ and equation
\eqref{eq:fTG}.
\end{proof}

If $\Delta$ is a simplicial complex, define
$$
F_{\Delta}(t) := \sum_{i} f_{i-1}(\Delta)t^i.
$$
If $\phi(q)$ is a function of $q$, define
$$
D_q \phi := \phi(q) - \phi(1).
$$
The point of $D_q$ is that it helps us convert from $B$ and homology
(the effect on which of taking skeleta is hard to describe) to $F$ and
$f$-vectors (the effect on which of taking skeleta is easy to
describe) in the following lemma.

\begin{lemma}\label{th:SDq}
If $\Delta \subseteq \Delta'$ are simplicial complexes, then
$$
S_{(\Delta,\Delta')}
 = (1+t^{-1})D_q S''_{(\Delta,\Delta')}
    + F_{\Delta} - F_{\Delta'}.
$$
\end{lemma}
\begin{proof}
By Lemma \ref{th:precursor},
$$
F_{\Delta}(t) - F_{\Delta'}(t)
  = S_{(\Delta,\Delta')}(t,1)
  = (1+t^{-1})S''_{(\Delta,\Delta')}(t,1) + B_{(\Delta,\Delta')}(t).
$$
Thus
$$
B_{(\Delta,\Delta')}(t) = 
  -(1+t^{-1})S''_{(\Delta,\Delta')}(t,1) 
  + F_{\Delta}(t) - F_{\Delta'}(t),
$$
which, when plugged back into Lemma \ref{th:precursor}, yields the desired result.
\end{proof}

\subsection{Skeleta}\label{su:skeleta}
Recall the {\em $s$-skeleton} of a simplicial complex $\Delta$ is
$$
\Delta\skel{s} := \{F \in \Delta\colon \dim F \leq s\}.
$$
Also recall that a simplicial complex is {\em pure} if all its facets have
the same dimension.
The {\em pure $s$-skeleton} of a simplicial complex $\Delta$ is
$$
\Delta\pskel{s} := \{F \in \Delta\colon F \subseteq G, G \in \Delta, \dim G = s\}.
$$
In other words, $\Delta\pskel{s}$ is the subcomplex of $\Delta$
consisting of the $s$-dimensional faces of $\Delta$, and all their
subfaces. (See \cite[Definition 2.8]{BjornerWachs1}.)
The results of the following lemma are easy exercises.

\begin{lemma}\label{th:del.con.skel}
If $\Delta$ is a simplicial complex and $e$ is a vertex of $\Delta$, then
\begin{enumerate}
\item  $(\Delta-e)\skel{s}   = \Delta\skel{s}-e$;
\item  $(\Delta/e)\skel{s-1} = \Delta\skel{s}/e$;
\item  $\Delta_s             = (\Delta\pskel{s})_s$;
\item  $(\Delta-e)_s         = (\Delta\pskel{s}-e)_s$; and
\item  $(\Delta/e)_{s-1}     = (\Delta\pskel{s}/e)_{s-1}$.
\end{enumerate}
\end{lemma}

\begin{lemma}\label{th:little.s.d-1}
If $\dim \Delta' \leq d-1$, then 
$$
\s''_{d-1}(\Delta,\Delta') \circeq \s''_{d-1}(\Delta\skel{d-1},{\Delta'}\skel{d-2}).
$$
\end{lemma}
\begin{proof}
Since $\Delta$ and $\Delta\skel{d-1}$ agree in dimensions $d-1$ and below,
$$
\s''_{d-1}(\Delta,\Delta') = \s''_{d-1}(\Delta\skel{d-1},\Delta').
$$
Next, replacing $\Delta'$ by ${\Delta'}\skel{d-2}$ in
$(\Delta\skel{d-1},\Delta')$ has the effect of adding
$(d-1)$-dimensional faces (in fact, all the $(d-1)$-dimensional faces
of $\Delta'$) to the simplicial pair, all of whose boundary faces are
still not present in the simplicial pair, since 
$\dim \Delta' \leq d-1$.
Thus
$$
\bdmap_{(\Delta^{(d-1)},\Delta'^{(d-2)});d-1} = 
  \bdmap_{(\Delta^{(d-1)},\Delta');d-1} \oplus 0
$$
(equivalently, the matrices representing the two boundary operators differ only in some additional zero columns); \cf\ proof of Lemma \ref{th:Ld}.  
It is then easy to check that, since $L''_{d-1} = \bdmap_{d-1}^*\bdmap_{d-1}$,
$$
L''_{d-1}(\Delta\skel{d-1},{\Delta'}\skel{d-2}) 
  = L''_{d-1}(\Delta\skel{d-1},\Delta') \oplus 0,
$$
and so
$$
s''_{d-1}(\Delta\skel{d-1},{\Delta'}\skel{d-2}) \circeq
 s''_{d-1}(\Delta\skel{d-1},\Delta').
$$
\end{proof}

\begin{corollary}\label{th:S.d-1}
If $\dim \Delta' \leq d-1$, then
$$
S''_{(\Delta,\Delta'),d-1} = S''_{(\Delta\skel{d-1},{\Delta'}\skel{d-2}),d-1}
$$
\end{corollary}

\begin{corollary}\label{th:S''.d-1}
If $\dim \Delta \leq d$ and $\dim \Delta' \leq d-1$, then
$$
S''_{(\Delta\skel{d-1},{\Delta'}\skel{d-2})} 
 = S''_{(\Delta,\Delta')} - S''_{(\Delta,\Delta'),d} t^{d+1}.
$$
\end{corollary}
\begin{proof}
Clearly, $(\Delta,\Delta')$ and
$(\Delta\skel{d-1},{\Delta'}\skel{d-2})$ agree in dimensions $d-2$ and
below.  Corollary \ref{th:S.d-1} thus ensures 
$
S''_{(\Delta\skel{d-1},{\Delta'}\skel{d-2})}
 = \sum_{i \leq d} S''_{(\Delta,\Delta'),i-1} t^i.
$
Then simply note, since $\dim \Delta \leq d$, that
$S''_{(\Delta,\Delta'),d}t^{d+1}$ is the only remaining term from
$S''_{(\Delta,\Delta')}$ not found in 
$S''_{(\Delta\skel{d-1},{\Delta'}\skel{d-2})}$.
\end{proof}

\subsection{The spectral recursion and skeleta}\label{su:spectral.skeleta}
\begin{lemma}\label{th:S.skeleton}
If $\dim \Delta \leq d$, $\dim \Delta' \leq d-1$, and $\Delta'
\subseteq \Delta$, then
\begin{multline*}
S_{(\Delta\skel{d-1},{\Delta'}\skel{d-2})}\\
 = S_{(\Delta,\Delta')} - (f_d(\Delta)+f_{d-1}(\Delta))t^{d+1}
 - (D_q S''_{(\Delta,\Delta'),d} - f_{d-1}(\Delta'))(t^{d+1} + t^d).
\end{multline*}
\end{lemma}
\begin{proof}
First 
use the definition of $D_q$ and
Corollary \ref{th:S''.d-1}
to get
\begin{equation}\label{eq:DqS''}
D_q S''_{(\Delta\skel{d-1},{\Delta'}\skel{d-2})}
 = D_q S''_{(\Delta,\Delta')} - (D_q S''_{(\Delta,\Delta'),d})t^{d+1}.
\end{equation}
Then apply Lemma \ref{th:SDq} (twice) and equation \eqref{eq:DqS''} to compute
\begin{align}
S_{(\Delta\skel{d-1},{\Delta'}\skel{d-2})}
 &=S_{(\Delta,\Delta')} 
   \begin{aligned}[t]
     &- (F_{\Delta} - F_{\Delta'})  
	  - (t^d+t^{d+1})D_q S''_{(\Delta,\Delta'),d}\\
     &+(F_{\Delta\skel{d-1}} - F_{{\Delta'}\skel{d-2}})
   \end{aligned} \nonumber \\
 &=S_{(\Delta,\Delta')} - (f_d(\Delta)t^{d+1} - f_{d-1}(\Delta')t^d)  
   - (t^d+t^{d+1})D_q S''_{(\Delta,\Delta'),d}. \label{eq:SftD}
\end{align}
The lemma now follows by adding the quantity
$(t^d+t^{d+1})f_{d-1}(\Delta')$ to the middle term of the right hand
side of equation \eqref{eq:SftD}, while subtracting it from the last
term.
\end{proof}

If $\Delta$ is a simplicial complex and $e$ is a vertex of $\Delta$, let
\begin{align*}
\calS_{\Delta,e} &:= S_\Delta - 
   (qS_{\Delta-e} + qtS_{\Delta/e} + (1-q)S_{(\Delta-e,\Delta/e)}), \\
\calS^d_{\Delta,e} &:= S''_{\Delta,d} - 
   (qS''_{\Delta-e,d} + qS''_{\Delta/e,d-1} 
    + (1-q)S''_{(\Delta-e,\Delta/e),d}),\ \text{and} \\
\Sf^d_{\Delta,e} &:= D_q\calS^d_{\Delta,e} + (1-q)f_{d-1}(\Delta/e).
\end{align*}
We have defined $\calS_{\Delta,e}$ precisely so that $\Delta$
satisfies the spectral recursion with respect to $e$ if and only if
$\calS_{\Delta,e}=0$, and we have defined $\calS^d_{\Delta,e}$ to be
the $d$-dimensional finer Laplacian version of $\calS_{\Delta,e}$.
The significance of $\Sf$ is made apparent by the next lemma, which is
the last key step to proving Theorem \ref{th:skeleta.preserve}.

\begin{lemma}\label{th:scriptS.skel}
If $\dim \Delta \leq d$ and $e$ is a vertex of $\Delta$, then 
$$
\calS_{\Delta\skel{d-1},e} = \calS_{\Delta,e} - 
   (t^d+t^{d+1})\Sf^d_{\Delta,e}.
$$
\end{lemma}
\begin{proof}
Since $\dim \Delta \leq d$, then $\dim \Delta-e \leq d$ and $\dim
\Delta/e \leq d-1$.  Therefore
\begin{multline*}
\calS_{\Delta\skel{d-1},e}\\
\begin{aligned}
 &= S_{\Delta\skel{d-1}} - qS_{(\Delta-e)\skel{d-1}} - qtS_{(\Delta/e)\skel{d-2}} 
      - (1-q)S_{((\Delta-e)\skel{d-1},(\Delta/e)\skel{d-2})}\\
 &= S_\Delta - qS_{\Delta-e} - qtS_{\Delta/e} - (1-q)S_{(\Delta-e,\Delta/e)}\\
 &\quad -f_d(\Delta)t^{d+1} + qf_d(\Delta-e)t^{d+1} + qtf_{d-1}(\Delta/e)t^d\\ 
 &\quad\quad + (1-q)(f_d(\Delta-e)+f_{d-1}(\Delta/e))t^{d+1}\\
 &\quad -   D_q S''_{\Delta,d}(t^{d+1}+t^d)
	 +  qD_q S''_{\Delta-e,d}(t^{d+1}+t^d)
	 + qtD_q S''_{\Delta/e,d-1}(t^d+t^{d-1})\\
 &\quad\quad +(1-q)(D_q S''_{(\Delta-e,\Delta/e),d} - f_{d-1}(\Delta/e))(t^{d+1}+t^d).
\end{aligned}
\end{multline*}
The first equation above is by the definition of $\calS$ and Lemma
\ref{th:del.con.skel}.  The second equation involves expanding each
term of the left-hand side by Lemma \ref{th:S.skeleton}, and then
regrouping like terms. Now, the second line and third lines of this
last expression add up to zero, by equation \eqref{eq:fTG}. The lemma
then follows from the definitions of $\calS$ and $\calS^d$.
\end{proof}

\begin{lemma}\label{th:Sf}
If $\dim \Delta \leq d$ and $e$ is a vertex of $\Delta$, then
$\calS_{\Delta,e}=0$ implies $\Sf^d_{\Delta,e}=0$.
\end{lemma}
\begin{proof}
It is easy to see that $\calS_{\Delta\skel{d-1},e}$ has no power of
$t$ higher than $\dim \Delta\skel{d-1} + 1 = d$.  But since
$\calS_{\Delta,e}=0$, Lemma \ref{th:scriptS.skel} implies that $0 =
\coeff{t^{d+1}}\calS_{\Delta\skel{d-1},e} = \Sf^d_{\Delta,e}$.  Here,
we are using the coefficient notation $\coeff{t^i}(\sum_j a_j t^j) := a_i$.
\end{proof}

\begin{lemma}\label{th:Sf.pureskel}
If $\Delta$ is a simplicial complex and $e$ is a vertex of $\Delta$, then 
$$
\Sf^d_{\Delta,e} = \Sf^d_{\Delta\pskel{d},e}.
$$
\end{lemma}
\begin{proof}
By expanding $\Sf^d_{\Delta,e}$ we need only show that
we may replace $\Delta$ by $\Delta\pskel{d}$ in each of
$S''_{\Delta,d}$,
$S''_{\Delta-e,d}$,
$S''_{\Delta/e,d-1}$,
$S''_{(\Delta-e,\Delta/e),d}$, and 
$f_{d-1}(\Delta/e)$.
But this follows from Lemma \ref{th:del.con.skel} and the definition of $S''$.
\end{proof}

\begin{theorem}\label{th:skeleta.preserve}
If $\dim \Delta \leq d$, and $e$ is a vertex of $\Delta$, then
$\Delta$ satisfies the spectral recursion with respect to $e$ iff 
$\Delta\skel{d-1}$ and $\Delta\pskel{d}$ do as well.
\end{theorem}
\begin{proof}
First assume $\Delta$ satisfies the spectral recursion with respect to $e$.  Then
$0=\calS_{\Delta,e}$.  By Lemma \ref{th:Sf}, then,
$\Sf^d_{\Delta,e}=0$.  And then by Lemma \ref{th:scriptS.skel},
$\calS_{\Delta\skel{d-1},e}=0$.  Furthermore,
$\Sf^d_{\Delta\pskel{d},e} = \Sf^d_{\Delta,e} = 0$,
by Lemma \ref{th:Sf.pureskel}.

Conversely, assume $\Delta\skel{d-1}$ and
$\Delta\pskel{d}$ satisfy the spectral recursion with respect to $e$.  By Lemma
\ref{th:Sf}, then $\Sf^d_{\Delta\pskel{d},e} = 0$.  And then by Lemmas
\ref{th:Sf.pureskel} and \ref{th:scriptS.skel},
$$
\calS_{\Delta,e} = \calS_{\Delta\skel{d-1},e} + (t^d + t^{d+1})\Sf^d_{\Delta,e}
 = \calS_{\Delta\skel{d-1},e} + (t^d + t^{d+1})\Sf^d_{\Delta\pskel{d},e} = 0.
$$
\end{proof}

\begin{corollary}\label{th:skeleta.preserve.all}
If $\dim \Delta \leq d$, then $\Delta$ satisfies the spectral recursion iff
$\Delta\skel{d-1}$ and $\Delta\pskel{d}$ do as well.
\end{corollary}

\subsection{Shifted complexes}\label{su:shift.finish}
Recall
a {\em $k$-set} is a set with $k$ elements, and a {\em $k$-family}
over ground set $E$ is a collection of $k$-subsets of $E$.  For a
$k$-set $F$, let $\bd F$ denote the $(k-1)$-family of all
$(k-1)$-subsets of $F$.  For a $k$-family $\K$, its {\em unsigned
boundary} $\bd \K$ is the $(k-1)$-family $\cup_{F \in \K} \bd F$.

If $F=\{f_1<\cdots<f_k\}$ and $G=\{g_1<\cdots<g_k\}$ are $k$-subsets
of integers, then $F \leq_P G$ under the {\em componentwise partial
order} if $f_p \leq g_p$ for all $p$.
A $k$-family $\K$ is {\em shifted} if $F \leq_P G$
and $G \in \K$ together imply that $F \in \K$.  A simplicial
complex 
$\Delta$ is {\em shifted} if $\Delta_i$ is shifted for every $i$.
The useful properties of shifted families in the following lemma are
easy to verify.
\begin{lemma}\label{th:shift.easy}
If $\K_1$ and $\K_2$ are shifted families, then so are $\bd \K_1$ and $\K_1 \cap \K_2$.
\end{lemma}

We say that $\Delta$ is a {\em near-cone with apex $1$} if
$\bd(\Delta-1) \subseteq \Delta/1$, where $\bd$ denotes the usual
unsigned boundary complex consisting of all faces that are not facets.
Equivalently, $\Delta$ is a near-cone with apex 1 if $F - v \disun 1
\in \Delta$ whenever $F \in \Delta$, $1 \not\in F$, and $v \in F$.
(See, \eg, \cite{BjornerKalai} for more on near-cones.)  We omit the
easy proofs of the following two lemmas.

\begin{lemma}
\label{th:inductive-shifted-definition}
Let $\Delta$ be a simplicial complex on $[n]$.  Then $\Delta$ is
shifted if and only if $\Delta$ is a near-cone with apex $1$, and both
$\Delta-1, \Delta/1$ are shifted with respect to the ordered vertex
set $[2,n]$.
\end{lemma}

\begin{lemma}\label{th:pure.shifted}
If $\Delta$ is a pure $d$-dimensional near-cone with apex $1$, then
$$
\Delta = (1 * (\Delta-1))\skel{d}
$$
\end{lemma}

\begin{theorem}\label{th:shifted}
If $\Delta$ is a shifted simplicial complex, then $\Delta$ satisfies
the spectral recursion, equation \eqref{eq:big}.
\end{theorem}
\begin{proof}
The proof is by induction on the dimension and number of vertices of
$\Delta$.  The base cases, when $\dim \Delta = 0$ or $\Delta$ has one
vertex (a special case of $\dim \Delta = 0$, anyway) are easy to
check.

Assume $\dim \Delta = d \geq 1$.  By induction, $\Delta\skel{d-1}$
satisfies the spectral recursion.  By Corollary \ref{th:skeleta.preserve.all}, it
remains to show that $\Delta\pskel{d}$ satisfies the spectral recursion as well.

To this end, first note that $\Delta_d$, the family of facets of
$\Delta\pskel{d}$, is shifted; then, 
by Lemma \ref{th:shift.easy} and reverse induction on dimension,
$\Delta\pskel{d}$ is shifted. 
By definition, $\Delta\pskel{d}$ is also pure, so Lemma
\ref{th:pure.shifted} implies
$$
\Delta\pskel{d} = (1 * (\Delta\pskel{d}-1))\skel{d}.
$$
Since $\Delta\pskel{d}$ is shifted, $\Delta\pskel{d}-1$ is also
shifted, with one less vertex, and so satisfies the spectral recursion, by
induction.  Thus $1*(\Delta\pskel{d}-1)$ also satisfies the spectral
recursion by Proposition \ref{th:v.ok} and Corollary
\ref{th:join.ok.all}.  Then Corollary \ref{th:skeleta.preserve.all}
guarantees that $\Delta\pskel{d}$ satisfies the spectral recursion.
\end{proof}

\section{Arbitrary shifted simplicial pairs}\label{se:arb.shifted}
Merris \cite{Merris} found a simple description of the Laplacian
spectrum of a shifted graph ($2$-family), in terms of the degree
sequence of the graph.  This was generalized in \cite{DuvalReiner} to
shifted families, by suitably generalizing the notion of degree
sequence.  In this section, we extend both the theorem, and the notion
of degree sequence, to shifted family {\em pairs} (Theorem
\ref{th:sdT}).  As in \cite{DuvalReiner}, the technique is to find
identical recursive formulas, similar to those in \cite{DuvalReiner},
for the Laplacian spectrum (Corollary \ref{th:spectra.shifted.family})
and the generalized degree sequence (Lemma \ref{th:d.shifted}), in
Subsections \ref{su:shift.laplace} and \ref{su:degree}, respectively.
The two threads are tied together with the proof of Theorem
\ref{th:sdT} in Subsection \ref{su:GM}.  Along the way, we rely upon
tools developed in Section \ref{se:shifted}.

Grone and Merris \cite{GroneMerris} conjectured that Merris'
description of the spectrum of a shifted graph becomes a majorization
inequality for an {\em arbitrary} graph.  This was also generalized
from graphs to families (though still not proved) in
\cite{DuvalReiner}.  In Subsection \ref{su:GM}, we also further
extend this conjecture from families to family pairs (Conjecture
\ref{th:GM}).

\subsection{Laplacians}\label{su:shift.laplace}
Recall the definition of family in Subsection \ref{su:shift.finish}.
If (for some $k$), $\K$ and $\K'$ are a $k$-family and $(k-1)$-family,
respectively, on the same ground set of vertices, then we will say
$(\K,\K')$ is a {\em family pair}, but we set $(\K,\K')=(\K,\K'')$
when $(\bd \K) \cap \K' = (\bd \K) \cap \K''$ (more formally, then, a
family pair is an equivalence class on ordered pairs of families).
We will say $(\K,\K')$ is a {\em shifted family pair} when $\K$ is
shifted and $(\K,\K')=(\K,\K'')$ for some $\K''$ that is shifted on
the same ordered ground set as $\K$.

Let $C(\K;\reals)$ denote the {\em oriented chains} of $k$-family
$\K$, \ie, the formal $\reals$-linear sums of oriented faces $[F]$
such that $F \in \K$.  If $(\K,\K')$ is a family pair, then the {\em
boundary operator} $\bdmap_{(\K,\K')}\colon C(\K;\reals) \rightarrow
C((\bd \K) \setm \K'; \reals)$ is defined as it is for simplicial
complexes, except that the sum is now restricted to faces in $(\bd \K)
\setm \K'$.  Equivalently,
$\bdmap_{(\K,\K')}=\bdmap_{(\cx{\K},\cx{\K'});k-1}$, when $\K$ is a
$k$-family and $\K'$ is a $(k-1)$-family.
As with simplicial complexes, the boundary operator has an adjoint
$\bdmap^*_{(\K,\K')}$, so 
the matrices representing $\bdmap$ and $\bdmap^*$ in the natural bases 
are transposes of one another.
\begin{definition}
The {\em Laplacian} of $(\K,\K')$ is the map
$L(\K,\K')\colon C(\K;\reals) \rightarrow C(\K;\reals)$ defined by
$$
L(\K,\K') := \bdmap_{(\K,\K')}^* \bdmap_{(\K,\K')}.
$$
It immediately follows that 
\begin{equation}\label{eq:alt}
L(\K,\K') = L''_{k-1}(\cx{\K},\cx{\K'}),
\end{equation}
where $\cx{\K}$ denote the
pure $(k-1)$-dimensional simplicial complex whose facets are the
members of $k$-family $\K$.
\end{definition}
It should be clear that $\bdmap_{(\K,\K')}$, and hence $L(\K,\K')$, is
well-defined on family pairs; that is, $\bdmap_{(\K,\K')} =
\bdmap_{(\K,\K'')}$ and $L(\K,\K') = L(\K,\K'')$, when
$(\K,\K')=(\K,\K'')$.  Of course, we may always specialize to a single
family by letting $\K'=\emptyset$.

Recall that $\Delta_i$ denotes the $(i+1)$-family of $i$-dimensional
faces of simplicial complex $\Delta$.
\begin{lemma}\label{th:Ld}
If $\dim \Delta' \leq d-1$, then 
$$
L''_d(\Delta,\Delta') = L(\Delta_d,\Delta'_{d-1}).
$$
\end{lemma}
\begin{proof}
The boundary maps $\bdmap_{(\Delta,\Delta');d}$ and
$\bdmap_{(\Delta_d,\Delta'_{d-1})}$ used to define
$L''_d(\Delta,\Delta')$ and $L(\Delta_d,\Delta'_{d-1})$, respectively,
both act on $C_d(\Delta;\reals)$.  By the definitions of $L$ and
$L''_d$, then, it will suffice to show that, for any $F \in \Delta_d$,
\begin{equation}\label{eq:missing.boundary}
\bdmap_{(\Delta,\Delta');d}[F]
=\bdmap_{(\Delta_d,\Delta'_{d-1})}[F].
\end{equation}
Now, the only difference between the left-hand and right-hand sides of
this equation is that the left-hand side is a sum restricted to faces
in the set difference $\Delta_{d-1} \setm \Delta'_{d-1}$, and the
right-hand side is a sum restricted to faces in $(\bd \Delta_d) \setm
\Delta'_{d-1}$.  Since $\Delta$ is a simplicial complex, $\bd \Delta_d
\subseteq \Delta_{d-1}$, so the only difference between the two sums
is provided by faces in $\Delta_{d-1} \setm (\bd \Delta_d)$.  But any
such face will not be in $\bd F$, the unsigned boundary of $F$, and
thus not appear in the expression for the signed boundary map, anyway.
(Equivalently, the matrices representing $\bdmap_{(\Delta,\Delta');d}$
and $\bdmap_{(\Delta_d,\Delta'_{d-1})}$ differ only in extra $0$ rows
indexed by $(d-1)$-dimensional faces of $\Delta$ not contained in any
$d$-dimensional face of $\Delta$, and these extra $0$ rows do not
affect $L=\bdmap^*\bdmap$.)  This establishes equation
\eqref{eq:missing.boundary}, and hence the lemma.  (\Cf\ the proof of
Lemma \ref{th:little.s.d-1}).
\end{proof}

Lemma \ref{th:Ld}, Proposition \ref{th:s.s''}, and equation
\eqref{eq:alt} allow us to go back and forth between families and
complexes.

\begin{lemma}\label{th:qt}
If $(\Gamma,\Gamma')$ is a simplicial pair, then
$$
qtS_{(\Gamma,\Gamma')} = S''_{(1*\Gamma,1*\Gamma')}.
$$
\end{lemma}
\begin{proof}
We compute $S_{(1*\Gamma,1*\Gamma')}$ in two different ways.  Since
$1*\Gamma$ and $1*\Gamma'$ are cones, $(1*\Gamma,1*\Gamma')$ has
trivial homology, so $B_{(1*\Gamma,1*\Gamma')}=0$.  Thus, by Lemma
\ref{th:precursor},
$$
S_{(1*\Gamma,1*\Gamma')} 
 = (1+t^{-1})S''_{(1*\Gamma,1*\Gamma')} + B_{(1*\Gamma,1*\Gamma')} 
 =t^{-1}(1+t)S''_{(1*\Gamma,1*\Gamma')}.
$$
On the other hand, by Corollary \ref{th:S.join},
$$
S_{(1*\Gamma,1*\Gamma')} = S_{1*(\Gamma,\Gamma')} = q(1+t)S_{(\Gamma,\Gamma')}.
$$
The lemma now follows immediately.
\end{proof}

Define $\s(\K,\K')$ to be the multiset of eigenvalues of $L(\K,\K')$,
arranged in weakly decreasing order.  When $\s(\K,\K')$ consists of
non-negative integers, it is a partition.  We will use the notation of
\cite{Macdonald} for partitions, except that we will denote the {\em
conjugate} or {\em transpose} of partition $\blambda$ by $\blambda^T$.
In particular, $1^m = (m)^T$ denotes the partition consisting of $m$
$1$'s.  Recall from Subsection \ref{su:fine.laplace} the definitions
of $\circeq$ and $\cup$ for multisets, which apply equally well to
partitions and weakly decreasing sequences.

Recall the definition of near-cone from subsection \ref{su:shift.finish}.
\begin{lemma}\label{th:spectra.shifted}
If $\Delta' \subseteq \Delta$ are pure near-cones with apex $1$, and
$\dim \Delta = d$ and $\dim \Delta' = d-1$, then, as partitions,
$$
\s''_d(\Delta,\Delta') \circeq  1^{f_{d-1}(\Delta/1) - f_{d-1}(\Delta'-1)} 
   + (\s''_{d}(\Delta-1,\Delta'-1) \cup \s''_{d-1}(\Delta/1,\Delta'/1)).
$$
\end{lemma}
\begin{proof}
Recall the coefficient notation $\coeff{t^i}(\sum_j a_j t^j) := a_i$.
First note
\begin{align}
\coeff{t^i}S''_{(\Gamma,\Gamma')} &= S''_{(\Gamma, \Gamma'),i-1}\label{eq:tS}\\
\coeff{t^i}B_{(\Gamma,\Gamma')} &= \tilde{\beta}_{i-1}(\Gamma, \Gamma')\label{eq:tB}
\end{align}
for any simplicial pair $(\Gamma, \Gamma')$ and for any $i$. Then,
by Lemmas \ref{th:S.d-1}, \ref{th:pure.shifted}, \ref{th:qt}, and equation \eqref{eq:tS},
$$
S''_{(\Delta,\Delta'),d} = q\coeff{t^d}S_{(\Delta-1,\Delta'-1)},
$$
so $\s''_d(\Delta,\Delta')$ has just as many non-zero parts as
there are terms in $q\coeff{t^d}S_{(\Delta-1,\Delta'-1)}$.
Lemma \ref{th:precursor} and equations \eqref{eq:tS} and \eqref{eq:tB} now imply
\begin{align*}
S''_{(\Delta,\Delta'),d}
 &= q\coeff{t^d}S_{(\Delta-1,\Delta'-1)}\\
 &= q(S''_{(\Delta-1,\Delta'-1),d-1} + S''_{(\Delta-1,\Delta'-1),d} +
        \tilde{\beta}_{d-1}(\Delta-1,\Delta'-1)),
\end{align*}
so the non-zero parts of $\s''_d(\Delta,\Delta')$ are given by 
adding $1$ to every element of the multiset union of three partitions:
$\s''_{d-1}(\Delta-1,\Delta'-1)$; $\s''_d(\Delta-1,\Delta'-1)$; and
the partition consisting of 
$\tilde{\beta}_{d-1}(\Delta-1,\Delta'-1)$ zeros.  This means
\begin{equation}\label{eq:spectra.shifted.0}
\s''_d(\Delta,\Delta') \circeq  1^m 
   + (\s''_{d-1}(\Delta-1,\Delta'-1) \cup \s''_{d}(\Delta-1,\Delta'-1)),
\end{equation}
where $m$ is the number of terms in
$q\coeff{t^d}S_{(\Delta-1,\Delta'-1)}$, since we established above that
$\s''_d(\Delta,\Delta')$ has $m$ non-zero parts.  But $\Delta' \subseteq
\Delta$ easily implies $\Delta'-1 \subseteq \Delta -1$, and so there
are
$$
m = f_{d-1}(\Delta-1) - f_{d-1}(\Delta'-1)
$$
terms in $q\coeff{t^d}S_{(\Delta-1,\Delta'-1)}$.

It is easy to verify that, since $\Delta$ and $\Delta'$ are pure 
near-cones (of dimensions $d$ and $d-1$, respectively) with apex $1$, 
\begin{align}
(\Delta-1)\skel{d-1}  &= \Delta/1;\ \text{and} \label{eq:spectra.shifted.noprime} \\
(\Delta'-1)\skel{d-2} &= \Delta'/1. \label{eq:spectra.shifted.prime}
\end{align}
From equation \eqref{eq:spectra.shifted.noprime}, we conclude
\begin{equation}\label{eq:spectra.shifted.f}
m = f_{d-1}(\Delta-1) - f_{d-1}(\Delta'-1) 
  = f_{d-1}(\Delta/1) - f_{d-1}(\Delta'-1).
\end{equation}
From equations \eqref{eq:spectra.shifted.noprime} and \eqref{eq:spectra.shifted.prime}, and
Lemma \ref{th:little.s.d-1}, we conclude
\begin{equation}\label{eq:spectra.shifted.last}
\s''_{d-1}(\Delta-1,\Delta'-1) \circeq \s''_{d-1}(\Delta/1,\Delta'/1).
\end{equation}
The lemma now follows from equations \eqref{eq:spectra.shifted.0},
\eqref{eq:spectra.shifted.f}, and \eqref{eq:spectra.shifted.last}.
\end{proof}

\begin{definition}
Let $\K$ be a $k$-family on ground set $E$, and $e \in E$. Then
the {\em deletion} of $\K$ with respect to $e$ is the $k$-family
$$
\K - e =\{F \in \K\colon e \not\in F\}
$$
on ground set $E-e$,
and the {\em contraction} of $\K$ with respect to $e$ 
is the $(k-1)$-family
$$
\K/e =\{F-e\colon F \in \K,\ e \in F\}
$$
on ground set $E-e$.
\end{definition}
The following identities are immediate:
$(\cx{\K}-e)_{k-1} = \K - e$,
$(\cx{\K}/e)_{k-2} = \K/e$, and
$\cx{\K}_{k-1} = \K$.

Define a $k$-family to be a {\em near-cone with apex $1$} when $\bd
(\K-1) \subseteq \K/1$.  It is an easy exercise to verify that $\K$ is
a near-cone iff $\cx{\K}$ is a near-cone.  Also, as with simplicial
complexes (Lemma \ref{th:inductive-shifted-definition}), $\K$ is
shifted iff $\K$ is a near-cone with apex $1$ such that $\K-1$ and
$\K/1$ are shifted.
The following corollary generalizes \cite[Lemma 5.3]{DuvalReiner}.
\begin{corollary}\label{th:spectra.shifted.family}
If $\K$ and $\K'$ are near-cone families with apex $1$ such that $\K'
\subseteq \bd\K$, then
$$
\s(\K,\K') \circeq 1^{\abs{\K/1} - \abs{\K'-1}}
                   + (\s(\K-1,\K'-1) \cup \s(\K/1,\K'/1)).
$$
\end{corollary}

\begin{proof}
Say $\K$ is a $k$-family, so $\K'$ is a $(k-1)$-family.  
Let $\Delta=\cx{\K}$ and $\Delta'=\cx{\K'}$.
From $\K' \subseteq \bd \K$, it follows that $\Delta' \subseteq \Delta$.
Then, by Lemmas \ref{th:Ld} and
\ref{th:spectra.shifted},
\begin{align*}
\s(\K,\K') 
 &= \s(\Delta_{k-1}, \Delta'_{k-2}) 
  = \s''_{k-1}(\Delta, \Delta') \\
 &\circeq 1^{f_{k-2}(\Delta/1) - f_{k-2}(\Delta'-1)}
	    + (\s''_{k-1}(\Delta-1,\Delta'-1) 
	     \cup \s''_{k-2}(\Delta/1,\Delta'/1)) \\
 &= 1^{f_{k-2}(\Delta/1) - f_{k-2}(\Delta'-1)} \\
	    &\quad + (\s((\Delta-1)_{k-1},(\Delta'-1)_{k-2})
	     \cup \s((\Delta/1)_{k-2},(\Delta'/1)_{k-3})) \\
 &= 1^{\abs{\K/1} - \abs{\K'-1}}
		    + (\s(\K-1,\K'-1) \cup \s(\K/1,\K'/1)). 
\end{align*}
\end{proof}

\subsection{Degree sequences}\label{su:degree}

\begin{notation}
We will write $F-\lambda$ to denote the set difference 
$F\backslash\{\lambda\}$, with the implicit assumption that $\lambda \in F$,
just as writing $F \disun \mu$ carries the implicit assumption that 
$\mu \not\in F$.  For instance, 
$\{F \in \K \colon F - \lambda \not\in \K'\}$ 
in the following definition is shorthand for
$\{F \in \K \colon \lambda \in F, F \backslash\{\lambda\} \not\in \K'\}$.
\end{notation}

\begin{definition}
Let $(\K,\K')$ be a family pair on ground set $E$.  Define the {\em
degree} of $\lambda$ in $(\K,\K')$ by
$$
d_\lambda(\K,\K') 
 :=\abs{\{F \in \K \colon F - \lambda \not\in \K'\}  }.
$$
It is easy to see that $d_\lambda$ is well-defined on family pairs;
that is, $d_\lambda(\K,\K')=d_\lambda(\K,\K'')$ when
$(\K,\K')=(\K,\K'')$.  The {\em degree sequence} $\dd=\dd(\K,\K')$ is
the partition whose parts are $\{d_\lambda\colon \lambda \in E\}$.
\end{definition}

In other words, to find the degree sequence of $(\K,\K')$, label all
the edges in the Hasse diagram of $\cx{\K}$ in the natural way, by the
vertex being added; then $d_\lambda$ counts the number of edges in the
Hasse diagram labelled $\lambda$, and connecting a face in $\K$ with a
face in $(\bd \K)\setm \K'$.
When $\K'=\emptyset$, then $\dd(\K)=\dd(\K,\emptyset)$ is the
generalized degree sequence of family $\K$ defined in \cite[Section
2]{DuvalReiner}.  It is also easy to see that $d_\lambda(\K) =
\abs{\K/\lambda}$.  When $\K$ is the set of edges of a graph, then
$\dd(\K)$ is the usual degree sequence of a graph.

\begin{lemma}\label{th:d.ordered}
If $(\K,\K')$ is a shifted family pair on $[1,n]$ and $1 \leq \lambda
< \mu \leq n$, then $d_\lambda(\K,\K') \geq d_\mu(\K,\K')$; \ie
$$\dd(\K,\K') = (d_1(\K,\K'),d_2(\K,\K'),\ldots,d_n(\K,\K')).$$ In
other words, the ordering of the degrees of the degree sequence of a
shifted family pair is given by the linear ordering of their vertices.
\end{lemma}

\begin{proof}
It will suffice to find an injection from $\{F \in \K\colon F - \mu
\not\in \K'\}$, a set whose cardinality equals $d_\mu(\K,\K')$, into $\{F
\in \K\colon F - \lambda \not\in \K'\}$, a set whose cardinality
equals $d_\lambda(\K,\K')$.  
It is easy to verify, using that $\K$ and $\K'$ are shifted, that such an 
injection $\phi$ is given by
$$
\phi(F) =
\begin{cases}
F			&\text{if $\lambda \in F$}\\
F - \mu \disun \lambda	&\text{if $\lambda \not\in F$}.
\end{cases}
$$
\end{proof}

The following lemma generalizes \cite[Lemma 5.2]{DuvalReiner}

\begin{lemma}\label{th:d.shifted}
If $\K$ and $\K'$ are shifted families on ground set $[1,n]$, and $\K'
\subseteq \bd\K$, then, as partitions,
$$
\dd(\K,\K')^T = 1^{\abs{\K/1} - \abs{\K'-1}}
  +(\dd(\K-1,\K'-1)^T \cup \dd(\K/1,\K'/1)^T).
$$ 
\end{lemma}
\begin{proof}
By standard partition arguments, 
this reduces to showing 
$$
\dd(\K,\K') = (\abs{\K/1} - \abs{\K'-1}) \cup
  (\dd(\K-1,\K'-1) + \dd(\K/1,\K'/1)),
$$ 
which 
is a direct consequence of the 
following two facts:
\begin{itemize}
\item $d_1(\K,\K') = \abs{\K/1} - \abs{\K'-1}$; and
\item if $\lambda>1$, then $d_\lambda(\K,\K') =  
  d_{\lambda}(\K-1,\K'-1) + d_{\lambda}(\K/1,\K'/1).$
\end{itemize}
The indexing on the second fact is indeed what is necessary, thanks to
Lemma \ref{th:d.ordered}, because $\K-1$, $\K'-1$, $\K/1$, and $\K'/1$
each have ground set $[2,n]$.
Each fact is an easy exercise, the first of which depends upon $\K$ being shifted.
\end{proof}

\subsection{A relative generalized Merris theorem}\label{su:GM}
Merris \cite[Theorem 2]{Merris} showed that when $\K$ is the
$2$-family of edges of a shifted graph, then $\s(\K) \circeq
\dd(\K)^T$.  This was generalized in \cite[Theorem 1.1]{DuvalReiner}
to allow $\K$ to be any shifted family.  The main result of this
section, below, further generalizes this to shifted family {\em
pairs}.  The proof is similar to that of \cite[Theorem
1.1]{DuvalReiner}.
\begin{theorem}\label{th:sdT}
If $(\K,\K')$ is a shifted family pair, then
$$
\s(\K,\K') \circeq \dd(\K,\K')^T
$$
\end{theorem}
\begin{proof}
By Lemmas \ref{th:shift.easy} and
\ref{th:inductive-shifted-definition}, $(\K-1,\K'-1) = (\K-1, (\K'-1)
\cap \bd (\K-1))$ and $(\K/1,\K'/1) = (\K/1, (\K'/1) \cap \bd (\K/1))$
are shifted family pairs.  Then the result is immediate from Corollary
\ref{th:spectra.shifted.family}, Lemmas
\ref{th:inductive-shifted-definition} and \ref{th:d.shifted}, and
induction on the number of vertices.
\end{proof}

Grone and Merris \cite[Conjecture 2]{GroneMerris} conjectured that
when $\K$ is the 2-family of edges of an {\em arbitrary} graph, then
the equality (modulo zeros) $\s(\K) \circeq \dd(\K)^T$ above becomes a
majorization inequality $\s(\K) \majby \dd(\K)^T$, \ie,
$\sum_{j=1}^k s_j \leq \sum_{j=1}^k d_j^T$ for all $k$, 
where $\s(\K)=(s_1,s_2,\ldots)$ and $\dd(\K)^T=(d_1^T,d_2^T,\ldots)$
are written as weakly decreasing sequences.  This majorization
inequality was also conjectured (but not proved) to hold when $\K$ is
any family, in \cite[Conjecture 1.2]{DuvalReiner}.  Based on no more
than a few examples, and that \cite[Theorem 1]{DuvalReiner}
successfully extends to pairs in Theorem \ref{th:sdT} above, we extend
this conjecture to family pairs as well.

\begin{conjecture}\label{th:GM}
If $(\K,\K')$ is a family pair, then
$$
\s(\K,\K') \majby \dd(\K,\K')^T.
$$
\end{conjecture}

Stephen \cite[Theorem 4.3.1]{Stephen} has shown that if the
Grone-Merris conjecture is true for all graphs, then Conjecture 5.9
holds for graph pairs ($\K$ is a 2-family and $\K'$ is a 1-family).

\begin{remark}\label{rm:fam.to.cx}
Theorem \ref{th:sdT} suffices to find the spectrum of a shifted
simplicial pair (that is, a simplicial pair $(\Delta,\Delta')$, where
$\Delta$ and $\Delta'$ are each shifted on the same ordered ground
set), not just a shifted family pair.  To see this, first note that by
Proposition \ref{th:s.s''}, finding $\s''_i(\Delta,\Delta')$ for all
$i$ determines the spectrum of the simplicial pair $(\Delta,\Delta')$.
Since $\s''_i$ depends only on $i$- and $(i-1)$-dimensional faces,
$\s''_i(\Delta,\Delta') = \s''_i(\Delta\skel{i},{\Delta'}\skel{i})$.
Finally, then,
$
\s''_i(\Delta,\Delta') 
 = \s''_i(\Delta\skel{i},{\Delta'}\skel{i})
 = \s''_i(\Delta\skel{i},{\Delta'}\skel{i-1})
 = \s(\Delta_i,\Delta'_{i-1})
$,
by Lemmas \ref{th:little.s.d-1} and \ref{th:Ld}.
\end{remark}

\section{Operations that preserve the spectral recursion}\label{se:preserve}
In this section, we see how the spectral recursion, equation
\eqref{eq:big}, and the spectrum polynomial behave with respect to
some natural operators on simplicial complexes.  Each operator has
significance for, or motivation from, matroids and/or shifted
complexes. Our main results are that the property of satisfying the
spectral recursion is preserved by disjoint union (Corollary
\ref{th:union.preserve.all}), Alexander duality (Corollary
\ref{th:spectral.Alexander.all}), and, with a slight modification
allowing order filters as well as simplicial complexes, two other dual
operators (Theorems \ref{th:spectral.dual} and
\ref{th:spectral.complement}).

\subsection{Duals}\label{su:dual}
The Tutte polynomial for matroids (see,
\eg, \cite{BrylawskiOxley}) whose recursion ($T_M=T_{M-e}+T_{M/e}$)
inspired and resembles the spectral recursion, is well-behaved with respect
to matroid duals ($T_M(x,y)=T_{M^*}(y,x)$), so it is natural to ask
what duality does to the spectrum polynomial and the spectral recursion.
There are three natural involutions on simplicial complexes that are
each appropriate generalizations of matroid duality.
How these involutions affect the
Laplacians of {\em families} has already been considered in
\cite[Section 4]{DuvalReiner}.  Recall that an {\em order filter}
$\Psi$ with vertices $V$ is a collection of subsets of $V$, closed
under taking supersets; that is, $F \in \Psi$ and $F \subseteq G
\subseteq V$ together imply $G \in \Psi$.
\begin{definition}
Let $\Delta$ be a simplicial complex (respectively, order filter) with
vertex set $V$.  The {\em dual} of $\Delta$ is the order filter (respectively,
simplicial complex)
$$
\Delta^*=\{V-F\colon F \in \Delta\}.
$$
The {\em complement} of $\Delta$ is the order filter (respectively,
simplicial complex)
$$
\Delta^c = \{F \subseteq V\colon F \not\in \Delta\}.
$$
The {\em Alexander dual} of $\Delta$ is the simplicial complex (respectively,
order filter)
$$
\Delta^\vee = \Delta^{*c} = \Delta^{c*}.
$$
\end{definition}
The Alexander dual has received attention lately in combinatorial
topology (see, \eg, \cite{BBLSW,BjornerButler}) and in combinatorial
commutative algebra (see, \eg, \cite{Bayer, BayerCharalambousPopescu,
BrunsHerzog, EagonReiner}).

It is easy to see that
$\Delta^{**}=\Delta^{cc}=\Delta^{\vee\vee}=\Delta$ for every
simplicial complex $\Delta$, and similarly for order filters.  If we
define an order filter $\Psi$ to be shifted when its every family
$\Psi_i$ of $i$-dimensional faces is shifted, then it is easy to see
that duality and complementation preserve being shifted, though with
the reverse vertex order.  Consequently, Alexander duality preserves
being shifted.

If $\Psi$ and $\Psi'$ are order filters on the same ground set of
vertices, we define the {\em order filter pair} $(\Psi,\Psi')$ to be
the simplicial pair $(\Psi'^c,\Psi^c)$, as defined in Section
\ref{se:laplace}.  (This means that, more formally, an order filter
pair is an equivalence class on ordered pairs of order filters.)  Thus
$(\Psi,\Psi')=(\Omega,\Omega')$ when the set differences $\Psi \setm
\Psi'$ and $\Omega \setm \Omega'$ are equal as subsets of the power set of
the ground set of vertices.  As with simplicial complexes, results and
definitions about order filter pairs $(\Psi,\Psi')$ may be specialized
to a single order filter, by letting $\Psi'=\emptyset$, the empty
order filter.

The definitions of deletion and contraction extend naturally to order
filters.  The deletion and contraction $\Psi-e$ and $\Psi/e$ of an
order filter $\Psi$ on vertex set $V$ are still order filters, though
on vertex set $V-e$.  
In contrast to simplicial complexes, $\Psi/e$ is not necessarily a subset
of $\Psi$ (though $\Psi-e \subseteq \Psi$, still), and $\Psi-e \subseteq \Psi/e$
(whereas, for simplicial complexes, $\Delta/e \subseteq \Delta -e$).

We now borrow a trick from \cite[Proposition 6]{KookReinerStanton}
(see also \cite[Proposition 4.2]{DuvalReiner}) to investigate how the
dual affects Laplacians and the spectral recursion.  Let
$(\Delta,\Delta')$ be a simplicial pair with
vertex set $[n]$; it is easy to specialize from pairs of duals to a
single dual, since the dual of the empty simplicial complex is again
empty, so $\Delta^* = (\Delta^*,\emptyset) = (\Delta^*,\emptyset^*)$.
Define $\phi_i(\Delta,\Delta')\colon C_i(\Delta,\Delta';\reals)
\rightarrow C_{n-i-2}(\Delta^*,\Delta'^*;\reals)$ to be the
$\reals$-linear isomorphism induced by
$$
\phi_i(\Delta,\Delta')\colon [F] \mapsto \sigma(F) [\overline{F}],
$$
where $\sigma(F)= (-1)^{\sum_{j \in F}j}$, and $\overline{F}=[n]-F$.

\begin{lemma}\label{th:fine.boundary.dual}
Let $(\Delta,\Delta')$ be a simplicial pair
with vertex set $[n]$, and let $\phi_j = \phi_j(\Delta,\Delta')$ for
any $j$.  Then
\begin{enumerate}
\item
$\phi_{i+1}^{-1} \bdmap_{(\Delta^*,\Delta'^*);n-i-2} \phi_i = - \bdmap^*_{(\Delta,\Delta');i+1}$, and
\item
$\phi_{i} \bdmap_{(\Delta,\Delta');i+1} \phi_{i+1}^{-1} = - \bdmap^*_{(\Delta^*,\Delta'^*);n-i-2}$.
\end{enumerate}
\end{lemma}
\begin{proof}
These are each a routine check of signs.
\end{proof}

\begin{corollary}\label{th:fine.L.dual}
Let $(\Delta,\Delta')$ be a simplicial pair
with vertex set $[n]$, and let $\phi_j = \phi_j(\Delta,\Delta')$ for
any $j$.  Then
$$
L_i(\Delta,\Delta') = \phi^{-1} L_{n-i-2}(\Delta^*,\Delta'^*)\phi.
$$
\end{corollary}

An immediate corollary is that, as first conjectured by V. Reiner
(personal communication),
\begin{equation}\label{eq:s.dual}
\s_i(\Delta,\Delta') = \s_{n-i-2}(\Delta^*,\Delta'^*),
\end{equation}
which translates into generating functions as
\begin{equation}\label{eq:S.tn.dual}
S_{(\Delta^*,\Delta'^*)}(t,q) = t^n S_{(\Delta,\Delta')}(t^{-1},q).
\end{equation}

We might hope that, if simplicial complex $\Delta$ satisfies the spectral
recursion with respect to a vertex $e$, then $\Delta^*$ would, too,
but this is not quite true.  
Routine calculations using equation \eqref{eq:S.tn.dual},
and duality identitites $(\Delta-e)^*=\Delta^*/e$ and $(\Delta/e)^*=\Delta^*-e$,
show that
\begin{equation}\label{eq:dual.spectral}
S_{\Delta^*}(t,q) 
   =qt S_{\Delta^*/e}(t,q) + 
    qS_{\Delta^*-e}(t,q) +
    (1-q) t S_{(\Delta^*/e,\Delta^*-e)}(t,q).
\end{equation}
We thus call
\begin{equation}\label{eq:spectral.filter}
S_{\Psi} = 
   qS_{\Psi-e}(t,q) +
   qt S_{\Psi/e}(t,q) + 
   (1-q) t S_{(\Psi/e,\Psi-e)}(t,q)
\end{equation}
{\em the spectral recursion for order filters}.  Theorem
\ref{th:spectral.complement} below provides further evidence that this
is the right formulation for order filters.  A unified approach to the
spectral recursions for simplicial complexes and order filters is to
develop a spectral recursion for simplicial complex pairs (which
includes simplicial complexes and order filters as special cases),
which is explored in \cite{Duval:RelLap}.

\begin{theorem}\label{th:spectral.dual}
If $\Delta$ is a simplicial complex and $e$ is an element of its
vertex set, then $\Delta$ satisfies the spectral recursion with respect to
$e$ iff $\Delta^*$ satisfies the spectral recursion for order filters,
equation \eqref{eq:spectral.filter}, with respect to $e$.
\end{theorem}
\begin{proof}
The forward implication follows from equation \eqref{eq:dual.spectral}
above.  The proof of the reverse implication is similar.
\end{proof}

The following proposition is a restatement of \cite[Corollary
4.7]{DuvalReiner}.

\begin{proposition}\label{th:m.Alex}
Let $\Delta$ be a simplicial complex with vertex set $[n]$.  If $\lambda \neq n$, then
$m_{\lambda}(L_i(\Delta)) = m_{\lambda}(L_{n-i-3}(\Delta^{\vee}))$.
\end{proposition}

The following corollary was first conjectured by V. Reiner (personal
communication).  

\begin{corollary}\label{th:complement}
If $\Delta$ is a simplicial complex with vertex set $[n]$, then
$\s_{i-1}(\Delta)$ and $\s_i(\Delta^c)$ agree, except for the
multiplicity of $n$.
\end{corollary}
\begin{proof}
By equation \eqref{eq:s.dual}, 
$
\s_i(\Delta^c) = \s_i(\Delta^{\vee*}) = \s_{n-i-2}(\Delta^{\vee})$,
so, if $\lambda \neq n$, then 
$$
m_{\lambda}(L_i(\Delta^c))
 = m_{\lambda}(L_{n-i-2}(\Delta^{\vee})
 = m_{\lambda}(L_{n-3-(n-i-2)}(\Delta))
 = m_{\lambda}(L_{i-1})
$$
by Proposition \ref{th:m.Alex}.
\end{proof}

The preceding proof is not as simple as it seems.  The proof of
Proposition \ref{th:m.Alex} in \cite[Corollary 4.7]{DuvalReiner} is
somewhat involved, and gets to the Alexander dual via the complement.
Especially in light of the simplicity of the statement of Corollary
\ref{th:complement}, we might hope it would have a more direct proof that does
not call upon the Alexander dual.

Corollary \ref{th:complement} translates into generating functions as
\begin{equation}\label{eq:S.complement}
S_{\Delta^c}(t,q)= tS_{\Delta}(t,q) + q^n A_{\Delta}(t),
\end{equation}
which we may rewrite as
\begin{equation}\label{eq:S.complement.alt}
S_{\Delta}(t,q) = t^{-1}S_{\Delta^c}(t,q) - q^n t^{-1}A_{\Delta}(t),
\end{equation}
where $A_{\Delta}(t)$ is a polynomial in $t$ that depends on $\Delta$.

\begin{theorem}\label{th:spectral.complement}
If $e$ is a vertex of simplicial complex $\Delta$, then $\Delta$
satisfies the spectral recursion with respect to $e$ iff $\Delta^c$
satisfies the spectral recursion for order filters, equation
\eqref{eq:spectral.filter}, with respect to $e$.
\end{theorem}
\begin{proof}
First assume $\Delta^c$ satisfies the spectral recursion for order
filters with respect to $e$.  Then, we may use equations
\eqref{eq:S.complement} and \eqref{eq:S.complement.alt},
and the complement identities 
$(\Delta-e)^c=\Delta^c-e$ and $(\Delta/e)^c=\Delta^c/e$,
to compute
\begin{equation}\label{eq:A.Delta}
S_{\Delta} = qS_{\Delta-e} + qtS_{\Delta/e} + (1-q)S_{(\Delta-e,\Delta/e)}
      +q^n t^{-1} (A_{\Delta-e} + A_{\Delta/e} - A_{\Delta}).
\end{equation}
By Lemma \ref{th:homology}, all simplicial complexes satisfy the spectral
recursion when $q=1$, so plugging $q=1$ into the above equation yields
$$
S_{\Delta}(1,t) = S_{\Delta}(1,t) + 1^n t^{-1} (A_{\Delta-e} + A_{\Delta/e} - A_{\Delta}).
$$
Therefore $A_{\Delta-e} + A_{\Delta/e} - A_{\Delta} = 0$, which, when
plugged back into equation \eqref{eq:A.Delta}, proves $\Delta$
satisfies the spectral recursion with respect to $e$.

The reverse implication is proved similarly.
\end{proof}

Theorems \ref{th:spectral.dual} and \ref{th:spectral.complement} together
imply the corresponding result for Alexander duality:
\begin{theorem}\label{th:spectral.Alexander}
If $e$ is a vertex of simplicial complex $\Delta$, then $\Delta$
satisfies the spectral recursion with respect to $e$ iff $\Delta^{\vee}$
satisfies the spectral recursion with respect to $e$.
\end{theorem}

\begin{corollary}\label{th:spectral.Alexander.all}
If $\Delta$ is a simplicial complex, then $\Delta$ satisfies the spectral
recursion iff $\Delta^{\vee}$ does as well.
\end{corollary}

\subsection{Union}\label{su:union}
The Tutte polynomial is well-behaved with respect to matroid direct
sum ($T_{M \oplus N} = T_M + T_N$), which corresponds to the union of
simplicial complexes with disjoint vertex sets ($\IN(M \oplus N) =
\IN(M) \cup \IN(N)$).  So it is natural to ask what disjoint union
does to the spectrum polynomial and the spectral recursion.

\begin{lemma}\label{th:S.union}
If $\Delta$ and $\Gamma$ are two non-empty simplicial 
complexes with disjoint vertex sets, then
$$
S_{\Delta \cup \Gamma} = S_{\Delta} + S_{\Gamma} + (t^{0} + 
t^{1})(q^{n+m} - 
(q^{n}+q^{m})) + t^{1}q^{0},
$$
where $\Delta$ and $\Gamma$ have $n=f_{0}(\Delta)$ and
$m=f_{0}(\Gamma)$ non-loop vertices, respectively.
\end{lemma}
\begin{proof}
For $i > 1$, it is clear that 
$L_{i-1}(\Delta \cup \Gamma) =   L_{i-1}(\Delta) \oplus 
L_{i-1}(\Gamma)$, since no $(i-1)$-dimensional face of $\Delta$ has 
any boundary in $\Gamma$, and vice versa.  Thus
$$
\s_{i-1}(\Delta \cup \Gamma) = \s_{i-1}(\Delta) \cup \s_{i-1}(\Gamma)
$$
for $i>1$.

The vertices of $\Delta$ and $\Gamma$ are disjoint, but they {\em 
share} the empty face in their boundary.  It is easy to see that 
$\s''_{0}(\Sigma) \circeq (f_{0}(\Sigma))$ for any simplicial complex 
$\Sigma$, so 
$\s''_{0}(\Delta \cup \Gamma) = (n+m)$, 
while 
$\s''_{0}(\Delta) \cup \s''_{0}(\Gamma) = (n,m)$.
Also, since 
$\dim L_{0}(\Delta \cup \Gamma) = f_{0}(\Delta \cup \Gamma) 
  = n+m = f_{0}(\Delta)+f_{0}(\Gamma) = \dim L_{0}(\Delta) + \dim 
  L_{0}(\Gamma)$, then
$\s_{0}(\Delta \cup \Gamma)$ and $\s_{0}(\Delta) \cup \s_{0}(\Gamma)$ 
have the same number of parts.  By Proposition \ref{th:s.s''}, it then 
follows that
$$
\s_{0}(\Delta \cup \Gamma) \cup (n,m)
 = \s_{0}(\Delta) \cup \s_{0}(\Gamma) \cup 
 (n+m,0).
$$
(In other words, to change $\s_{0}(\Delta) \cup \s_{0}(\Gamma)$ into
$\s_{0}(\Delta \cup \Gamma)$, replace $(n,m)$ 
in $\s_{0}(\Delta) \cup \s_{0}(\Gamma)$ by 
$(n+m,0)$ in $\s_{0}(\Delta \cup \Gamma)$.)
Similarly, since $\Delta \cup \Gamma$, $\Delta$, and $\Gamma$ each 
have exactly one empty face, $\s_{-1}(\Delta \cup \Gamma)$ has one 
element, and $\s_{-1}(\Delta) \cup \s_{-1}(\Gamma)$ has two elements, 
and so
$$
\s_{-1}(\Delta \cup \Gamma) = (n+m),
$$
while
$$
\s_{-1}(\Delta) \cup \s_{-1}(\Gamma) = (n,m).
$$
The lemma now follows immediately.
\end{proof}    

We continue to assume $\Delta$ and $\Gamma$ are non-empty simplicial
complexes with disjoint vertex sets, and that 
$\Gamma$ has 
$m$ non-loop vertices.
By arguments similar to those in the proof of Lemma \ref{th:S.union},
$$
S_{(\Gamma,\emptyset)} = 
 S_{\Gamma} - (t^{0} +t^{1})q^{m} + t^{1}q^{0},
$$
and so
\begin{equation}\label{eq:rel.union}
S_{(\Delta \cup \Gamma,\Delta')} = S_{(\Delta,\Delta')} + S_{\Gamma}
  -(t^{0}+t^{1})q^{m} + t^{1}q^{0}.
\end{equation}
    
\begin{theorem}\label{th:union.preserve}
If $\Delta$ satisfies the spectral recursion with respect to $e$, and 
$\Gamma$ is any simplicial complex whose vertex set is disjoint from 
the vertex set of $\Delta$, then $\Delta \cup \Gamma$ satisfies the spectral 
recursion with respect to $e$.    
\end{theorem}
\begin{proof}
If $\Gamma=\emptyset$, then the theorem is trivially true.  
Otherwise, 
it is a routine calculation with
Lemma \ref{th:S.union} and equation \eqref{eq:rel.union}.
\end{proof}

\begin{corollary}\label{th:union.preserve.all}
If $\Delta$ and $\Gamma$ each satisfy the spectral recursion, then so does
their disjoint union $\Delta \cup \Gamma$.
\end{corollary}

The following example shows that the {\em arbitrary} union of two
simplicial complexes satisfying the spectral recursion does not itself
necessarily satisfy the spectral recursion, even if both complexes are pure.

\begin{example}
Let $\Delta$ be the pure $1$-dimensional simplicial complex on vertex
set $\{a,b,c,d,e\}$ with facets $\{ab,ac,ad,ae,bc,bd\}$.  (We omit
brackets and commas from each face for clarity.)  Let $\Gamma$ be the
pure $1$-dimensional simplicial complex on the same vertex set 
with facets $\{ab,ac,ad,ae,de\}$.  Now, $\Delta$ is shifted with
vertices ordered $a<b<c<d<e$, and $\Gamma$ is shifted with vertices
ordered $a<d<e<b<c$, so each satisfies the spectral recursion.

On the other hand, we can easily show $\Delta \cup \Gamma$ does not
satisfy the spectral recursion with respect to vertex $d$.  First check
directly that $\Delta \cup \Gamma$ is not Laplacian integral.  (Note
that $\Delta \cup \Gamma$ is the $1$-dimensional skeleton of the cone
over Example \ref{ex:universal}.)  Next, since $(\Delta \cup
\Gamma)-d$ and $(\Delta \cup \Gamma)/d$ are each isomorphic to shifted
complexes (with different vertex orders), they are each Laplacian integral.
It is also easy to directly verify that $((\Delta \cup \Gamma)-d,
(\Delta \cup \Gamma)/d)$ is Laplacian integral as well.  Thus, the
right-hand side of the spectral recursion in this instance has all integer
exponents, but the left-hand side does not.
\end{example}

\section{Acknowledgements}
I am grateful to Vic Reiner, Eric Babson, Woong Kook, and Tamon
Stephen for many helpful comments and suggestions, to Vic Reiner also
for introducing me to this problem, and to several referees for their
insightful advice.

\newcommand{\journalname}[1]{\textit{#1}}
\newcommand{\booktitle}[1]{\textit{#1}}

\end{document}